\documentclass{amsart}

\usepackage{amsmath,amssymb,graphicx,subfigure}
\newtheorem{thm}{Theorem}[section]
\newtheorem{prop}{Proposition}
\newtheorem{lem}{Lemma}
\newtheorem{cor}{Corollary}
\newtheorem{defn}{Definition}

\newcommand{\bi}{\begin{itemize}}
\newcommand{\ei}{\end{itemize}}
\newcommand{\be}{\begin{enumerate}}
\newcommand{\ee}{\end{enumerate}}
\newcommand{\bc}{\begin{center}}
\newcommand{\ec}{\end{center}}
\newcommand{\bt}{\begin{tabular}}
\newcommand{\et}{\end{tabular}}

\newcommand{\BS}{BS$(1,2)$}
\newcommand{\bs}{Baumslag-Solitar}
\newcommand{\cg}{Cayley graph}
\newcommand{\gset}{generating set}
\newcommand{\cf}{context-free}
\newcommand{\texp}{$t$-exponent}
\newcommand{\lan}{language}
\newcommand{\nf}{normal form}
\newcommand{\NF}{\mathcal{NF}}
\newcommand{\mc}{\mathcal}
\newcommand{\geo}{geodesic}

\newcommand{\ra}{\rightarrow}
\newcommand{\pda}{pushdown automaton}
\newcommand{\fsa}{finite state automaton}
\newcommand{\Z}{\mathbb Z}
\newcommand{\hnn}{HNN-extension}

\begin{document}
%%%%%%%%%%%%%%%%

\title {A context-free and a 1-counter geodesic language for a
Baumslag-Solitar~group}

\author {Murray Elder} \address{School~of
Mathematics~and~Statistics, University of St~Andrews, North~Haugh,
St~Andrews, Fife, KY16 9SS, Scotland } \email{murray@mcs.st-and.ac.uk}
\thanks{Supported by EPSRC grant GR/S53503/01 }

\begin{abstract}
We give a \lan\ of unique \geo\ \nf s for the \bs\ group \BS\ that is
context-free and 1-counter. We discuss the classes of \cf, 1-counter and
counter languages, and explain how they are inter-related.
\end{abstract}


  \date{\today} \maketitle

\section{Introduction}
In this article we give a simple combinatorial description of a \lan\
of \nf s for the solvable \bs\ group \BS\ with the standard \gset,
such that each \nf\ word is \geo, each group element has a unique \nf\
representative, and the \lan\ is accepted by a (partially blind)
1-counter automaton.  It follows that the language is \cf.

Several authors have studied \geo\ \lan s for the (solvable) \bs\
groups, including Brazil \cite{Brazil}, Collins, Edjvet and Gill
\cite{CollinsEG}, Freden and McCann \cite{Freden}, Groves
\cite{Groves}, Miller \cite{Miller}, and the author and Hermiller
\cite{EHerm}. It is well known that \bs\ groups are asynchronously
automatic but not automatic \cite{Epstein}, and the asynchronous \lan\
is not geodesic.  Groves proved that no geodesic \lan\ of \nf s for a
solvable \bs\ group with standard \gset\ can be regular \cite{Groves},
so we could say that context-free or 1-counter is the next best thing.

Collins, Edjvet and Gill proved that the growth function (the formal
power series where the $n$-th coeficient is the number of elements
having a geodesic representative of length $n$) of a solvable \bs\
group is rational \cite{CollinsEG}, and Freden and McCann have studied
growth functions for the non-solvable case \cite{Freden}.
%% cf. algebraic

If $G$ is a group with generating set $\mc G$, we say two words $u,v$
are equal in the group, or $u=_Gv$, if they represent the same group
element. We say $u$ and $v$ are identical if the are equal in the free
monoid, that is, they are equal in $\mc G^*$.

\begin{defn}[$G$-automaton]\label{defn:EE} 
Let $G$ be a group and $\Sigma$ a finite set.  A (non-deterministic)
 {\em $G$-automaton} $A_G$ over $\Sigma$ is a finite directed graph
 with a distinguished {\em start vertex} $q_0$, some distinguished
 {\em accept vertices}, and with edges labeled by $(\Sigma^{\pm
 1}\cup\{\epsilon\})\times G$.  If $p$ is a path in $A_G$, the element
 of $(\Sigma^{\pm 1})$ which is the first component of the label of
 $p$ is denoted by $w(p)$, and the element of $G$ which is the second
 component of the label of $p$ is denoted $g(p)$. If $p$ is the empty
 path, $g(p)$ is the identity element of $G$ and $w(p)$ is the empty
 word.  $A_G$ is said to {\em accept} a word $w\in (\Sigma^{\pm 1})$
 if there is a path $p$ from the start vertex to some accept vertex
 such that $w(p)=w$ and $g(p)=_G 1$.
 \end{defn}

\begin{defn}[Finite state automaton; Regular]
If $G$ is the trivial group, then $A_G$ is a (non-deterministic)
{\em finite state automaton}. A language is {\em regular} if it is the set
of strings accepted by a finite state automaton.
\end{defn}

\begin{defn}[Counter; 1-counter]
A language is {\em $k$-counter} if it is accepted by some
$\Z^k$-automaton. We call the generators of $\Z^k$ and their inverses
{\em counters}. A language is {\em counter} if it is $k$-counter for
some $k\geq 1$.
\end{defn}

\begin{figure}[ht!]
  \bc
     \includegraphics[height=3cm]{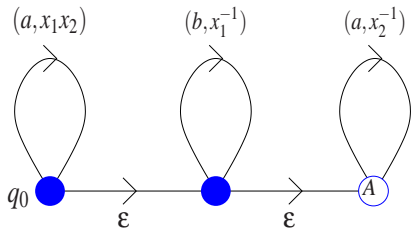}
  \ec
  \caption{A counter automaton accepting $a^nb^na^n$.}
  \label{fig:counterNotCF}
\end{figure}

For example, the language $\{a^nb^na^n \; | \; n\in \mathbb N\}$ is
accepted by the $\Z^2$-automaton in Figure \ref{fig:counterNotCF}, with
alphabet $a,b$ and counters $x_1,x_2$.

In the case of $\Z$-automata, we assume that the generator is $1$ and
the binary operation is addition, and we may insist without loss of
generality each transition changes the counter by either $0,1$ or
$-1$. We can do this by adding states and transitions to the automaton
appropriately.  That is, if some edge changes the counter by $k\neq
0,\pm 1$ then divide the edge into $|k|$ edges using more states.  The
symbols $+,-$ indicate a change of $1,-1$ respectively on a
transition.

%Note that with this definition, a $\Z$-automaton cannot ``see'' the
%value of the counter until it reaches an accept state, so is not
%equivalent to a \pda\ with one stack symbol, which can determine if
%the stack is empty at any time.  In Section \ref{sect:prelim} we show
%that this definition of 1-counter implies context-free. We also show
%that these 1-counter languages are not closed under concatenation.
%Recall that a {\em regular language} is the set of strings accepted by
%a {\em finite state automaton} (see \cite{Sipser}).

\begin{defn}[Pushdown automaton; Context-free]
A {\em pushdown automaton} is a 6-tuple $(Q,\Sigma,\Gamma, \tau, q_0,
A)$ where $Q,\Sigma, \Gamma$ and $A$ are all finite sets, and \be
\item $Q$ is the set of states,
\item $\Sigma$ is the input alphabet together with the empty word $\epsilon$,
\item $\Gamma$ is the stack alphabet together with $\epsilon$ (the
empty symbol),
\item $\tau$ is the transition function, 
\item $q_0$ is the start state,
\item $A\subseteq Q$ is the set of  accept states.\ee

The transition function takes as input a state and an input letter,
and outputs a state and a stack instruction of the form $\gamma \ra
\beta$, which means pop $\gamma$ from the top of the stack then push
$\beta$ on the top of the stack. Note that $\epsilon \ra \gamma$ means
push $\gamma$ onto the stack, $\gamma \ra \epsilon$ means pop $\gamma$
off the stack, and $\epsilon \ra \epsilon$ means do nothing (and in
this case will be omitted).

A word is accepted by the automaton if there is a sequence of
transitions starting from the state $q_0$ with an empty stack, pushing
and popping stack symbols, to an accept state. Note that you can
always push new symbols onto the stack, but you can only pop if the
correct symbol is on top of the stack.

A language is {\em context-free} if it is the language of some \pda.
\end{defn}

\begin{figure}[ht!]
  \bc
     \includegraphics[height=3cm]{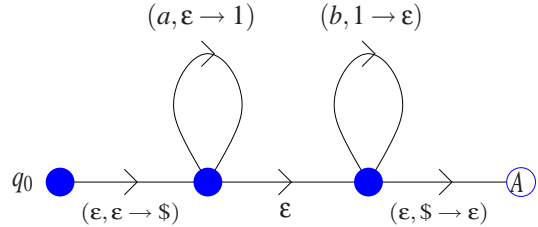}
  \ec
  \caption{Pushdown automaton accepting $a^nb^n$.}
  \label{fig:PDAanbn}
\end{figure}

As an example, the language $\{a^nb^n \; | \; n\in \mathbb N\}$ is
accepted by the \pda\ in Figure \ref{fig:PDAanbn} with alphabet $a,b$
and stack symbols $\$,1$, and this language is not regular
\cite{HU},\cite{Sipser}.

Note that our definition of counter automata is not equivalent to a
pushdown automata with a stack (with one type of token) for each
counter, since in our definition, we cannot test the value of the
counter until we are done reading the input. For this reason, these
automata are sometimes referred to as ``partially blind'' or
vision-impaired counter automata, since the cannot ``see'' whether the
counter is non-zero except at the end.

\begin{defn}[\bs\ group]
The group with presentation \\ $\langle a,t \; | \;
tat^{-1}=a^p\rangle$ is the {\em solvable \bs\ group} BS$(1,p)$, for
\\ $p\in \Z, p\geq 2$.
\end{defn}

\begin{figure}[ht!]
  \bc
 \includegraphics[width=13cm]{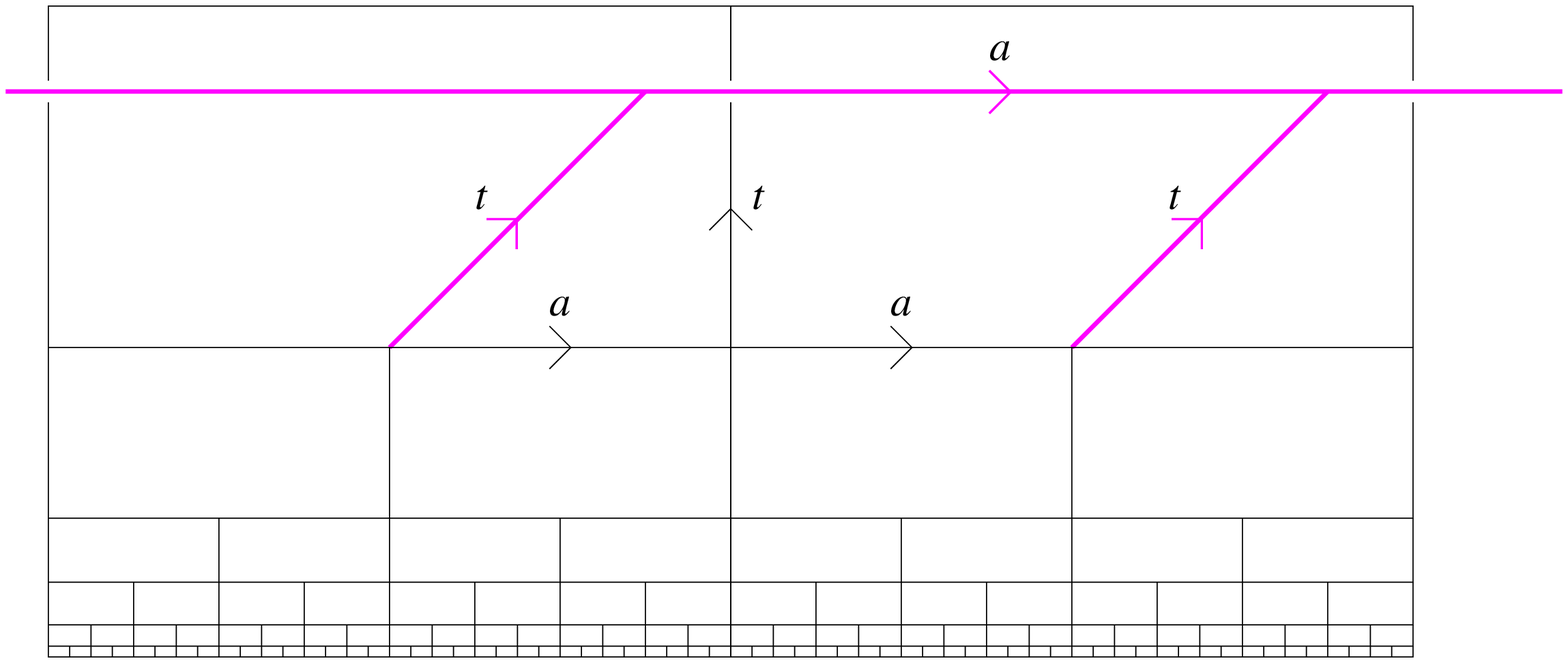}
  \ec
  \caption{Part of the \cg\ for \BS. }
  \label{fig:BS12}
\end{figure}

\noindent
In this article we will consider the group \BS. Let $\mc
G=\{a,a^{-1},t,t^{-1}\}$ be the inverse closed \gset\ for \BS.  We
give a picture of part of the \cg\ for \BS\ in Figure
\ref{fig:BS12}. From the side the \cg\ looks like a binary tree. See
\cite{Epstein} for a detailed description of the \cg.

The paper is organised as follows. In Sections \ref{sect:prelim} and
\ref{sect:counternotcf} we examine the various definitions of formal
languages presented above, and establish their relative intersections
and inclusions, which we illustrate in Figure \ref{fig:sets}. In
particular we prove that 1-counter languages as defined are \cf.  In
Section \ref{sect:nf} we define a normal form language for \BS\ and
prove that each normal form word is geodesic, and the language of
normal form words bijects to the set of group elements.  In Section
\ref{sect:mainthm} we prove that this normal form language is
1-counter, which implies it is context-free. Then in the last section
we show that the language of all geodesics for \BS\ is not counter.

\section{1-counter languages}\label{sect:prelim}

\begin{lem}\label{lem:1countercf}
Every 1-counter language is context-free.
\end{lem}

\noindent
\textit{Proof.} Let $L$ be a 1-counter language accepted by a
1-counter machine $M$. We will construct a (non-deterministic) \pda\ $N$
that accepts the language $L$, with stack symbols $\$_+,\$_-$ and
$1$. Let $M_+$ be a copy of $M$ obtained by replacing transitions
$(a,+)$ by $(a,\epsilon\ra 1)$ and $(a,-)$ by $(a,1\ra \epsilon)$, and
let $M_-$ be a copy of $M$ obtained by replacing transitions $(a',-)$
by $(a',\epsilon \ra 1)$ and $(a',+)$ by $(a',1\ra \epsilon)$.

$N$ is constructed from these two automata $M_+$ and $M_-$ as follows.
The states of $N$ consist of two distinct states $q_+,q_-$ for each
state $q$ of $M$, plus a new start state $s_0$ and a new single accept
state $p$.  There is a transition labelled $(\epsilon, \epsilon \ra
\$_+)$ from $s_0$ to the former start state $(q_0)_+$ in $M_+$. For each
$q_+$ in $M_+$ there is a transition labelled $(\epsilon, \$_+ \ra
\$_- )$ from $q_+$ to the corresponding state $q_-$ in $M_-$, and a
transition labelled $(\epsilon, \$_- \ra \$_+ )$ from $q_-$ to $q_+$
in $M_+$.

Finally for every accept state $q$ in $M$ there is a transition
labelled $(\epsilon, \$_+\ra \epsilon)$ from $q_+$ in $M_+$ to the
single accept state $p$, and $(\epsilon, \$_-\ra \epsilon)$ from $q_-$
in $M_-$ to the single accept state $p$.

This new machine works by starting with an empty stack and pushing
$\$_+ $ on the bottom. Then if the old machine increments the counter,
the new machine adds $1$ to the stack. From then on if the counter
value never dips below zero, the new machine will stay in the $M_+$
states. However if there is ever a ``pop 1'' but the symbol on the
stack is $\$_+ $, pass over to $M_-$. Then the height of the stack now
represents the negative value of the counter, you stay in this side
until the value of the counter comes back to zero, in which case you
can switch.

It follows that the language of $N$ is precisely the language of the
1-counter machine $L$.  $\Box$

%% Mitrana proves that L(F_2) =CFL

\begin{lem}\label{lem:not1counter2}
The language of strings of the form $a^mb^ma^nb^n$ is both counter and
context-free but not 1-counter.
\end{lem}

\begin{figure}[ht!]
  \bc
\bt{c}
\includegraphics[height=3cm]{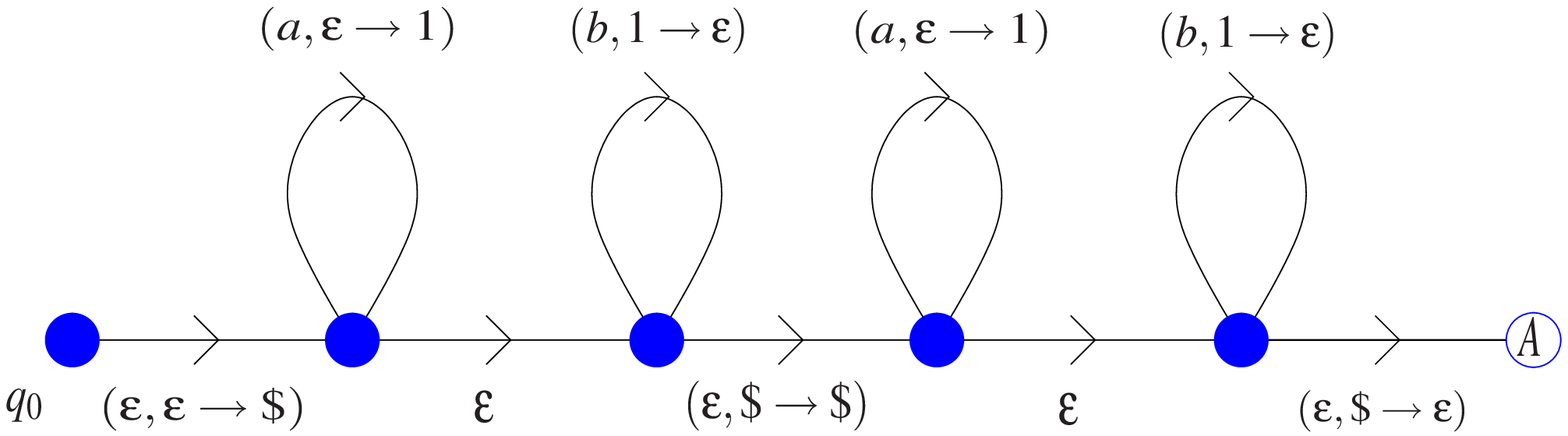}\\
\\
\includegraphics[height=3cm]{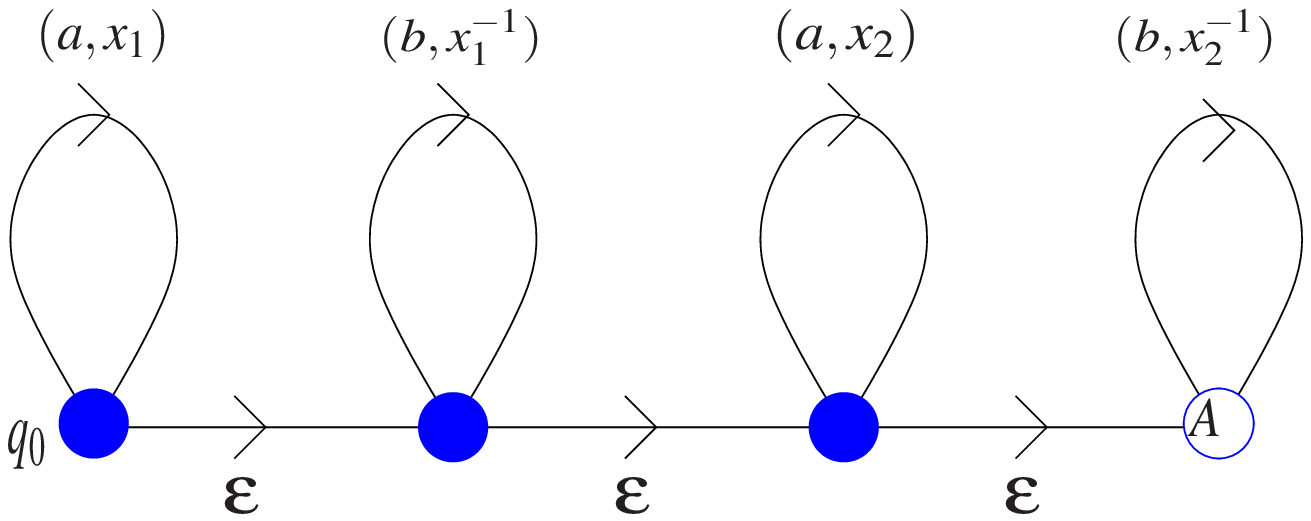}
\et
\ec
  \caption{Pushdown automaton and 2-counter machine accepting $a^mb^ma^nb^n$}
  \label{fig:4counter}
\end{figure}

\textit{Proof.} The \pda\ and the $\Z^2$-automaton in Figure
\ref{fig:4counter}
both accept this language, so it is \cf\ and counter.

Suppose by way  of contradiction that the language  is 1-counter, and
 let $M$ be a 1-counter machine for it with $p$ states. Assume without
 loss  of generality  that  each transition  changes  the counter  by
 either $0, -1$ or $1$.
 
Define $a_1=a^{p^2}, b_1=b^{p^2}, a_2=a^{p^2},b_2=b^{p^2}$, and
consider the word $s=a_1b_1a_2b_2$
%a^{p^2}b^{p^2}a^{p^2}b^{p^2}$ 
which belongs to the language.

Consider the prefix $a_1=a^{p^2}$. Since this prefix is longer than the
number of states, it must visit some state twice, so
$a_1=x_0y_0z_0$ where $y_0$ represents a loop of length at most
$p$.

If going around $y_0$ causes a net change of zero in the value of the
counter, then going around it twice would give a new word that is
accepted by $M$, but not of the form $a^mb^ma^nb^n$. So assume the net
change is $k_0$ with $|k_0|\geq 1$.

Let $s_1=x_0z_0$ which has length at least $p^2-p$, so must go around
a loop in $M$. So $s_1=x_1y_1z_1$ with $y_1$ a loop of length at most
$p$. Again, if the net change in the counter going around $y$ is zero
then we can go around $y_1$ twice and have a word accepted by $M$ that
is not in the language.

If the net change is $k_1$ of the opposite sign to $k_0$ then there is
a word that goes $|k_1|$ times around the loop $y_0$ then $|k_0|$
times around $y_1$, which keeps the final value of the counter at
zero, so is accepted by $M$, but since we are pumping the $a^{p^2}$
prefix of $s$ we have a word that is not in the language.

Thus $y_1$ changes the counter by $k_1$ with $|k_1|\geq 1$ and having the
same sign of $k_0$. Let $s_2=x_1z_1$ with length at least $p^2-2p$.

Iteratively we can write $s_i=x_iy_iz_i$ with $y_i$ a loop which
changes the value of the counter by an amount $k_i$ of the same sign
as $k_0$, until there are no loops left in $x_iz_i$, which does not
happen until at least $p$ iterations (since $s_i$ has length at least
$p^2-ip$).

Since $x_iz_i$ has no loops, it has length at most $p-1$. So it
changes the value of the counter by at most $l_1$ where
$|l_1|<p$. Whereas, the sum of the $|y_i|$ changes the value of the
counter by at least $p$ since each one contributes at least $1$ to the
sum.

Now repeat this analysis for the subwords $b_1,a_2$ and $b_2$.

If all the loops in each subword change the counter by the same sign,
then we have a contradiction, since the net change of all the loops is
greater than $4p$ whereas the net change of the four remaining
$x_iz_i$ segments is less than $4p$, so they cannot cancel each other.

Thus at least two subwords have loops of opposite signs.  If the loops
in $a_1$ have the same sign as the loops in $a_2$ and $b_2$, then the
loops in $b_1$ must have the opposite sign. So suppose that some loop
in $b_1$ changes the counter by $k$, and some loop in $a_2$ changes
the counter by $l$ of the opposite sign to $k$. Then pumping the first
loop by $|l|$ and the second by $|k|$ gives a word that is accepted by $M$
and not in the language.

Otherwise if the loops in $a_1$ have the opposite sign to the loops in
either $a_2$ or $b_2$, then take a loop in $a_1$ which changes the
counter by $k$ and a loop in $a_2$ or $b_2$ that changes the counter
by $l$ of the opposite sign to $k$. Then pumping the first loop by
$|l|$ and the second by $|k|$ gives a word that is accepted by $M$ and
not in the language.  $\Box$

\begin{cor}
1-counter languages are not closed under concatenation or intersection.
\end{cor}

\noindent
\textit{Proof.} The language $C=\{a^nb^n \; | \; n \in \mathbb N \}$
is 1-counter but $CC$ is not 1-counter by the previous lemma (Lemma
\ref{lem:not1counter2}).  

The languages $D=\{a^nb^nc^m \; | \; m,n \in \mathbb N \}$ and
$E=\{a^mb^nc^n \; | \; m,n \in \mathbb N \}$ are 1-counter, but $D\cap
E=\{a^nb^nc^n \; | \; n \in \mathbb N \}$ is not context-free
\cite{HU},\cite{Sipser} so by Lemma \ref{lem:1countercf} is not
1-counter. $\Box$

However, we have
\begin{lem}[Closure properties of $k$-counter \lan s]
\label{lem:closure1counter} 
If $C,C'$ are $k$-counter for $k\geq 1$ and $L$ is regular, then
$C\cup C'$, $C\cap L$, $CL$ and $LC$ are all $k$-counter.
\end{lem}

\noindent
\textit{Proof.}  Let $M,M'$ be $k$-counter automata for $C,C'$, with
start states $q_0,q_0'$, states $S,S'$, and accept states $A,A'$,
respectively.  Then construct a $k$-counter automaton accepting $C\cup
C'$ with a new start state $p_0$ joined to $q_0,q_0'$ by two epsilon
transitions.

Let $N$ be a \fsa\ for $L$ with states $T$, start state $p_0$ and
accept states $B$.  Construct a $k$-counter automaton accepting $C\cap
L$ having states $S\times T$, start state $(q_0,p_0)$, such that
$(q,p)$ is an accept state if $q\in A,p\in B$ (they are both accept
states), and if there are transitions from $q$ to $r$ in $M$ labelled
by $(a,g)$ and $p$ to $s$ in $N$ labelled by $a$ where $g\in \Z^k$,
then there is a transition from $(q,p)$ to $(r,s)$ labelled $(a,g)$.

Construct a $k$-counter automaton accepting $CL$ with start state
$q_0$ and accept states $B$ by adding an epsilon transition from each
accept state of $M$ to $p_0$.

Construct a $k$-counter automaton accepting $LC$ with start state
$p_0$ and accept states $A$ by adding an epsilon transition from each
accept state of $N$ to $q_0$.  $\Box$

Iterating the union operation a finite number of times gives
\begin{cor}\label{cor:union}
The union of a finite number of $k$-counter languages is $k$-counter.
\end{cor}

\section{Context-free and not counter}\label{sect:counternotcf}

The language $\{a^nb^na^n\; | \; n\in \mathbb N\}$ accepted by the
$\Z^2$-automaton in Figure \ref{fig:counterNotCF} is not \cf\ by standard
results \cite{HU},\cite{Sipser}.  In this section we show that
conversely, there is a language that is context-free but not counter.

Consider a string of letters $a,b,c$. We say a string contains a {\em
square} if it has a subword of the form $ww$.  An interesting result
from combinatorics is that one can write out a square-free word in
$a,b,c$ of arbitrary length. This is due to Thue and Morse and
described in \cite{Lothaire} (Chapter 2).  In particular we have
\begin{prop}[Thue-Morse]\label{prop:squarefree} Define a homomorphism
$f$ on $\{a,b,c\}$ by $f(a)=abc,f(b)=ac$ and $f(c)=b$.  Then for
any $i\in \mathbb N$, $f^i(a)$ is square-free.  \end{prop}

For example, to compute  $f^3(a)$ we have \\
$a\ra abc\ra abcacb\ra abcacbabcbac$. 

 In order to show that a language is not counter we make use of the
following lemma.
%than the ``interchange lemma'' given in \cite{Mitrana}.

\begin{lem}[Swapping Lemma] \label{lem:swap}
If $L$ is counter then there is a constant $s>0$, the ``swapping
length'', such that if $w\in L$ with length at least $2s+1$ then $w$
can be divided into four pieces $w=uxyz$ such that $|uxy|\leq 2s+1$,
$|x|,|y|>0$ and $uyxz\in L$.
\end{lem}

\noindent
\textit{Proof.}  Let $s$ be the number of states in the counter
automaton, and let $p$ be a path in the $Z^k$-automaton such that
$w(p)=w$. If $p$ visits each state at most twice then it cannot have
length more than $2s$, so $p$ visits some state at least three times.
Let $u$ be the first part of $w(p)$ until it hits this state, then $x$
a non-trivial loop back to this state the second time, $y$ a loop back
a third time, and $z$ the rest of $w$. So $w(p)=uxyz$ ends at an
accept state, and the second component of $p$ equals
$g(u)g(x)g(y)g(z)=_{\Z^k} 1$.  Switching the orders of $x$ and $y$,
the path $uyxz$ still takes you to the same accept state, and
$g(uyxz)=_{\Z^k} 1$ since all elements of $\Z^k$ commute, so $uyxz\in
L$.  $\Box$

Note its similarity to the pumping lemmas for regular and context-free
languages \cite{HU},\cite{Sipser}.  This lemma is only of any use if
your word $w$ has no squares, otherwise you can just swap the square
and get the same word (that is $x=y$).

\begin{thm}\label{thm:cfnotcounter}
There is a language that is context-free but not counter.
\end{thm}

\noindent
\textit{Proof.} Consider the language of all strings in $a,b,c$ of the
form $ww^R$, where $w^R$ is word obtained by reversing $w$. It is well
known that this is a \cf\ language \cite{HU},\cite{Sipser}, since it
is accepted by a \pda\ which uses the stack to store the first half of
the word, then checks the last half of the word matches.

Suppose by way of contradiction that this language is counter, with
swapping length $p$ as in Lemma \ref{lem:swap}.  Let $w$ be a
square-free word from Proposition \ref{prop:squarefree} of length at
least $2p+1$. Then $ww^R$ can be split into four subwords $u,x,y,z$
such that $uxy$ falls in the first $w$ prefix. Since $w$ has no
squares and $x,y$ are adjacent words then it must be that $x\neq
y$. But $uyxz$ will fail to be in the language because the second part
will not be the reverse of the first part.  $\Box$

\begin{figure}[ht!]
  \bc
  \includegraphics[width=13.5cm]{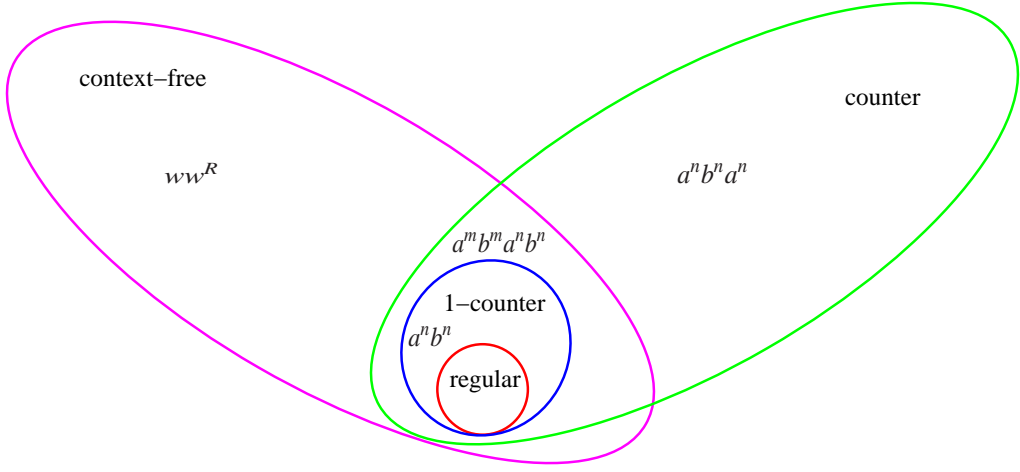}
  \ec
  \caption{Intersections of the formal languages}
  \label{fig:sets}
\end{figure}

In Figure \ref{fig:sets} we have a diagram of sets of regular,
1-counter, \cf\ and counter languages, and by the above results we
have shown the given inclusions.

The fact that there are counter languages that are not \cf\ and vice
versa can be observed by considering word problems for various groups.
The word problem for a group $G$ with generating set $\mc G$ is the
set $WP(G)=\{w\in \mc G^*: \overline w=1\}$ of all words in the
generating set that evaluate to the identity element. By work of
Muller and Schupp \cite{MS}, the word problem for the group $\Z^2$ is
not a \cf\ language, whereas the word problem of the free group on two
(or more) generators is \cf.  Elston and Ostheimer \cite{EO} proved
that a group has a deterministic counter word problem (with a so-called
inverse property) if and only if it is virtually abelian, so the word
problem for $Z^2$ is counter. To see why $WP(F_2)$ is not counter,
consider a Thue-Morse word made up of an arbitrary number of subwords
$(aaa),(aba),(ab^{-1}a)$, followed by its ``reverse'' in the subwords
$(a^{-1}a^{-1}a^{-1}),(a^{-1}b^{-1}a^{-1}),(a^{-1}ba^{-1})$. This word
is in the word problem, but applying the Swapping lemma (Lemma
\ref{lem:swap}) gives a word that is non-trivial.

%follows that the word problem for the group $\Z^2$ is a counter
%language but is not a \cf\ language, whereas the word problem of the
%free group on two (or more) generators is \cf\ but not counter
%\footnote[1]{In fact their results are much stronger than this. Muller
%and Schupp prove that a group has a context-free word problem if and
%only if it is virtually free, and Elston and Ostheimer prove that a
%group has a word problem that is accepted by a deterministic counter
%automaton with a so-called weak inverse property if and only if it is
%virtually abelian. Herbst \cite{Herbst} shows that a group has a
%1-counter word problem if and only if it is virtually cyclic. Finally,
%a group has a regular word problem if and only if it is finite.}.  

The first examples of languages that are counter but not context-free
were given by Mitrana and Stiebe in \cite{Mitrana}. Mitrana and Stiebe
give the following lemma, which they call the ``interchange lemma'',
which they use to show that the language of palindromes, and the
language $\{a^ib^i\;| i\geq 0\}^*$, are not counter. We include it here
for completeness, and to show how it differs from the Swapping Lemma
above.

\begin{lem}[Interchange Lemma \cite{Mitrana}]
If $L$ is the language of a $G$-automaton where $G$ is an abelian
group, then there is a constant $p$ such that for any word $x\in L$ of
length at least $p$, and for any given subdivision of $x$ into
subwords $v_1w_1v_2w_2 \ldots w_pv_{p+1}$ with $|w_i|\geq 1$, there
are some $r,s$ such that the word obtained from $x$ by interchanging
$w_r$ and $w_s$ is in $L$.
\end{lem}

\section{The \nf\ \lan}\label{sect:nf}

Recall that \BS\ $=\langle a,t \; | \; tat^{-1}=a^2\rangle$ with the
(standard) inverse closed \gset\ $\mc G=\{a,a^{-1},t,t^{-1}\}$.  We
wish to describe geodesic words with respect to this generating set.

\begin{defn}[$E,N,P,X$]
A word is of the form $E$ if it is $a^i$. A word is of the form $N$ if
it has no $t$ letters and at least one $t^{-1}$ letter. A word is of
the form $P$ if it has no $t^{-1}$ letters and at least one $t$
letter.

A word is of the form $X$ if it is the concatenation of a $P$ word of
$t$-exponent $k$, followed by an $N$ word of $t$-exponent $(-k)$. That
is, an $X$ word is a word of type $PN$ with zero \texp.
\end{defn}

\begin{figure}[ht!]
  \bc
               \includegraphics[width=13.5cm]{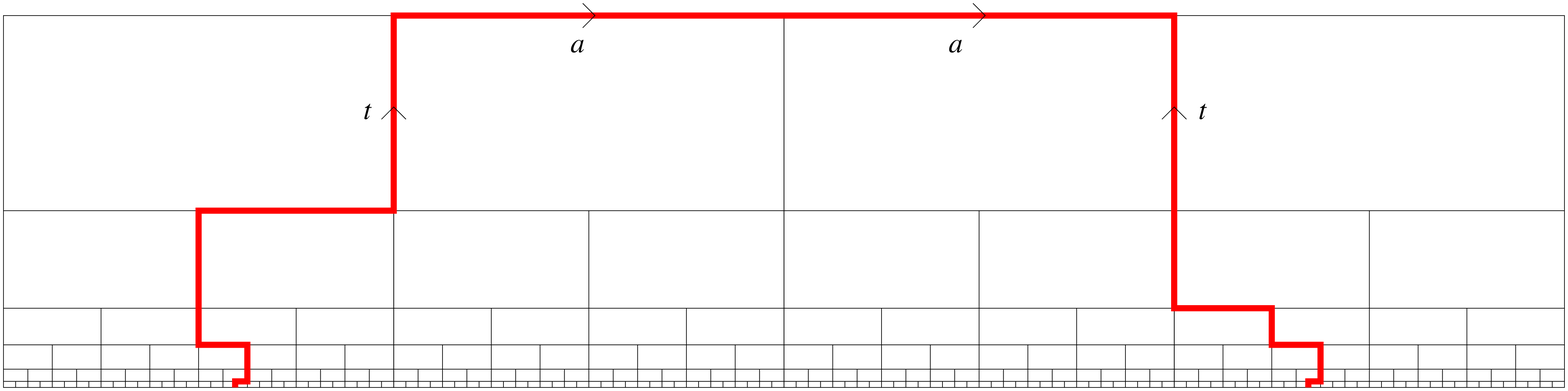}
  \ec
  \caption{An $X$ word}
  \label{fig:mesa}
\end{figure}

Benson Farb called words of type $X$ ``mesas'', since drawing an $X$
word in the \cg\ resembles this land formation. See Figure
\ref{fig:mesa}.

While the following fact is well known, we include an elementary proof
of it here for completeness.
\begin{lem}[Commutation]\label{lem:commutation}
If $u$ has zero \texp\ then $au=ua$ and $a^{-1}u=ua^{-1}$.
\end{lem}

\noindent
\textit{Proof.} If $u$ is type $X$ then $u=_{BS}a^i$ so
$au=a^{i+1}=ua$.

If $u$ is type $NP$ then let $u=vw$ where $v$ is type $N$ with \texp\
$-k$ and $w$ is type $P$ (so has \texp\ $k$). Each time we push $a^i$
past a $t^{-1}$ it becomes $a^{2i}$ since $at^{-1}=t^{-1}a^2$. Then
$au=avw=va^{2^k}w$. Each time we push $a^{2^i}$ past a $t$ it becomes
$a^{2^{i-1}}$ since $a^2t=ta$. So $au=avw=va^{2^k}w=vwa=ua$. Finally
if $u$ is any other form, first replace each occurrence of
$ta^it^{-1}$ in $u$ by $a^{2i}$.  Then $u$ becomes a word of type $NP$
with zero \texp. We can pass $a$ through this word as in the previous
case, and then put $u$ back in its original form and we are done.  $
\Box $

\begin{lem}[Miller \cite{Miller}]
 \label{lem:Miller} Every geodesic word in $\mc G
^*$ is a subword of a word of type $NPN$ or $PNP$.  \end{lem} See
Lemma $1$ of \cite{Groves} for a proof. We can use this lemma to
describe a subset of geodesic words that represent every group
element.

Define a type $NP_{\leq}$ word to be a word of type $NP$ with
non-positive \texp\ sum, and type $NP_>$ to be a word of type $NP$
with positive \texp\ sum.  

\begin{lem}[Ten types]\label{lem:TenTypes} Every
element of \BS\ has a geodesic representative in $\mc G^*$ that is one
of  ten types:\\ $E,X,N,XN,NP_{\leq}, XNP$ having \texp\ $\leq
0$, or \\ $P,PX,NP_>,NPX$ having \texp\ $>0$, such that no more than
three $a$ or $a^{-1}$ letters can occur in succession in the
geodesic.\end{lem}

Hermiller and the author used a similar characterisation in our work
on minimal almost convexity \cite{EHerm}.

\noindent
\textit{Proof.}  Every group element can be represented by some
geodesic word in $\mc G^*$.  If a geodesic word has no $t^{\pm 1}$
letters then it is type $E$. Otherwise by Lemma \ref{lem:Miller} it is
a word of type $N,P,NP,PN,NPN$ or $PNP$.

If the geodesic is type $NP$ then it either has non-positive
$t$-exponent sum, so is type $NP_{\leq}$, or positive $t$-exponent
sum, so is type $NP_>$.

If the geodesic is type $PN$ then it either has zero $t$-exponent sum,
so is type $X$, negative $t$-exponent sum, so is type $XN$, or
positive $t$-exponent sum, so is type $PX$.

Suppose the geodesic is a word $w$ of type $NPN$. If $w$ has positive
$t$-exponent sum it is type $NPX$. If $w$ has zero $t$-exponent sum,
then write it as $ux$ where $u$ is type $NP$ with zero $t$-exponent
sum and $x$ is type $X$. By Lemma \ref{lem:commutation} $w=_{BS}xu$
which has the same length and is type $XNP$.  If $w$ has negative
$t$-exponent sum, then
$w=a^{\epsilon_1}t^{-1}uta^{\epsilon_2}txt^{-1}a^{\epsilon_2}t^{-1}v$
where $u$ is type $E$ or $NP$ with zero $t$-exponent sum, $x$ is type
$E$ or $X$, $v$ is type $E$ or $N$, and $\epsilon_i\in \mathbb Z$.
Then by Lemma \ref{lem:commutation}
\begin{eqnarray*}
w  =_{BS}  a^{\epsilon_1+\epsilon_2+\epsilon_3}(t^{-1}ut)(txt^{-1})t^{-1}v \\
   =_{BS}  a^{\epsilon_1+\epsilon_2+\epsilon_3}(txt^{-1})(t^{-1}ut)t^{-1}v
\end{eqnarray*}
which is not geodesic since we can cancel $tt^{-1}$ at the end.

Finally, suppose the geodesic is a word $w$ of type $PNP$. If $w$ has
negative or zero $t$-exponent sum it is type $XNP$.  If $w$ has
positive $t$-exponent sum, then
$w=a^{\epsilon_1}txt^{-1}a^{\epsilon_2}t^{-1}uta^{\epsilon_2}tv$ where
$x$ is type $E$ or $X$, $u$ is type $E$ or $NP$ with zero $t$-exponent
sum, $v$ is type $E$ or $P$, and $\epsilon_i\in \mathbb Z$.  Then by
Lemma \ref{lem:commutation}
\begin{eqnarray*}
w  =_{BS}  a^{\epsilon_1+\epsilon_2+\epsilon_3}(txt^{-1})(t^{-1}ut)tv \\
  =_{BS}    a^{\epsilon_1+\epsilon_2+\epsilon_3}(t^{-1}ut)(txt^{-1})tv 
\end{eqnarray*}
which is not geodesic since we can cancel $t^{-1}t$ at the end.

The additional condition that no more than three $a$'s are allowed in
succession is obtained by observing that $a^6=_{BS}t^3at^{-1}$ so any
power of $a$ greater than five is not geodesic, and since
$a^4=ta^2t^{-1}$ and $a^5=ta^2t^{-1}a=ata^2t^{-1}$ we choose to
replace $a$-exponents of $4$ or $5$ by subwords of the same length.  An
identical argument eliminates powers of $a^{-1}$ greater than three.
$\Box$

\begin{defn}[Run]
An {\em $N$-run} is a word of the form
$$a^{\epsilon_k}t^{-1}a^{\epsilon_{k-1}}t^{-1}\ldots
t^{-1}a^{\epsilon_1}t^{-1}a^{\epsilon_0}.$$ A {\em $P$-run} is a word
of the form $$a^{\epsilon_0}ta^{\epsilon_1}t\ldots
ta^{\epsilon_{k-1}}ta^{\epsilon_k}.$$ We can write a run in shorthand
by just writing the $a$-exponents.
 For example,
$a^2t^{-1}at^{-1}a^0t^{-1}at^{-1}a^{-1}$ can be written as $2101(-1)$.

We call the $a$-exponents {\em entries} of the run. A run is {\em
non-trivial} if it has at least one non-zero entry. Note that a run
that has at least one $t$ or $t^{-1}$ letter will have at least two
entries, since by definition a run starts and ends with a power of $a$
(possibly $a^0$).

We say a geodesic has at most one non-trivial run if it can be
expressed as the concatenation of geodesic $N$- or $P$-runs such that
at most one factor is non-trivial. For example, the word
$t^2a^2t^{-1}at^{-2}$ can be written as $(t^2)(a^2t^{-1}at^{-2})$, so
has at most one run.  \end{defn}

Drawing the $N$-run represented by $2101(-1)$ in the \cg\ we start to
see what behaviour is allowed in a geodesic. For instance, the
sub-runs $1(-1)$ and $(-1)1$ are not allowed since
$$at^{-1}a^{-1}\ra t^{-1}a \;\;\;\;\;\;\; a^{-1}t^{-1}a\ra
t^{-1}a^{-1}.$$

Also, if the $N$-run $2101(-1)$ were preceded by a $t^{-1}$ then we
would have $t^{-1}a^2$ which can be written as $at^{-1}$.  In fact,
the only time you could ever see an entry that is not $0,1$ or $-1$ is
at the {\em start} of an $N$-run, or the {\em end} of a $P$-run.
\begin{figure}[ht!]
  \bc
               \includegraphics[width=9cm]{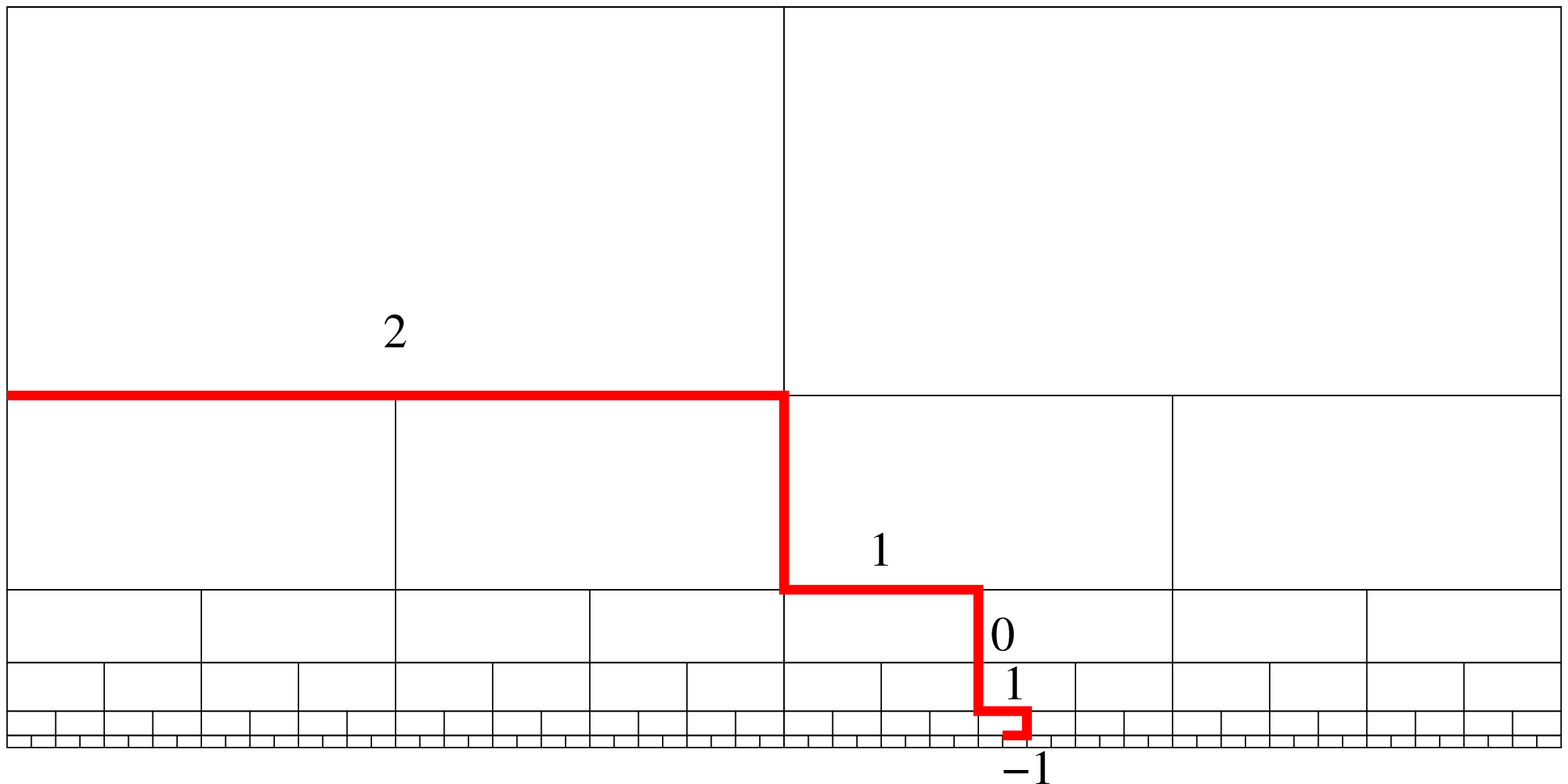}
  \ec
  \caption{The $N$-run $2101(-1)$}
  \label{fig:ExampleRun}
\end{figure}

\begin{lem}[No $|i|>6$] \label{lem:No6} 
If a run represents a geodesic word and has an entry $i$ that is not
one of $1,0$ and $(-1)$, then $i$ must be one of
$2,3,4,5,(-2)(-3),(-4),(-5)$ and occurs at the {\em start} of an
$N$-run or the {\em end} of a $P$-run.
\end{lem}

\noindent
\textit{Proof.}  If $i\geq 6$ occurs at any point in a run then
 $a^6\ra ta^3t^{-1}$ so the run is not geodesic.

For $N$-runs, if $i\geq 2$ occurs after the start of the run then
$t^{-1}a^2\ra at^{-1}$ so the run is not geodesic.  If $i\leq -2$
occurs after the start of the run then $t^{-1}a^{-2}\ra a^{-1}t^{-1}$
so the run is not geodesic.

For $P$-runs, if $i\geq 2$ occurs before the end of the run then
$a^2t\ra ta$ so the run is not geodesic. If $i\leq -2$ occurs before
the end of the run then $a^{-2}t\ra ta^{-1}$ so the run is not
geodesic.
$\Box$

\begin{lem}[No consecutive $1(-1),(-1)1$] \label{lem:No-11} 
A geodesic run cannot contain $1(-1)$ or $(-1)1$.
\end{lem}

\begin{figure}[ht!]
\bc
\bt{ccc}
\bt{c}
\includegraphics[width=4cm]{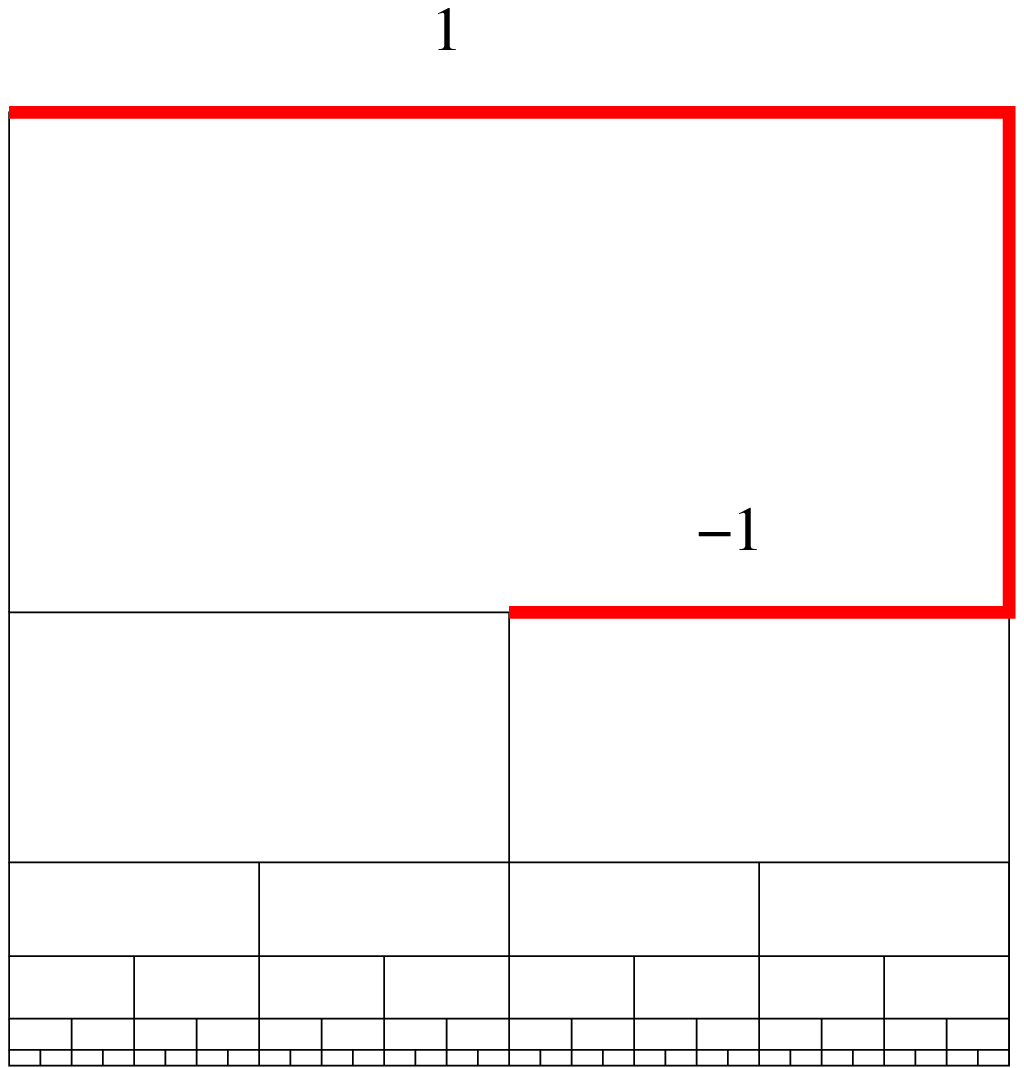}
\et& &
\bt{c}
\includegraphics[width=4cm]{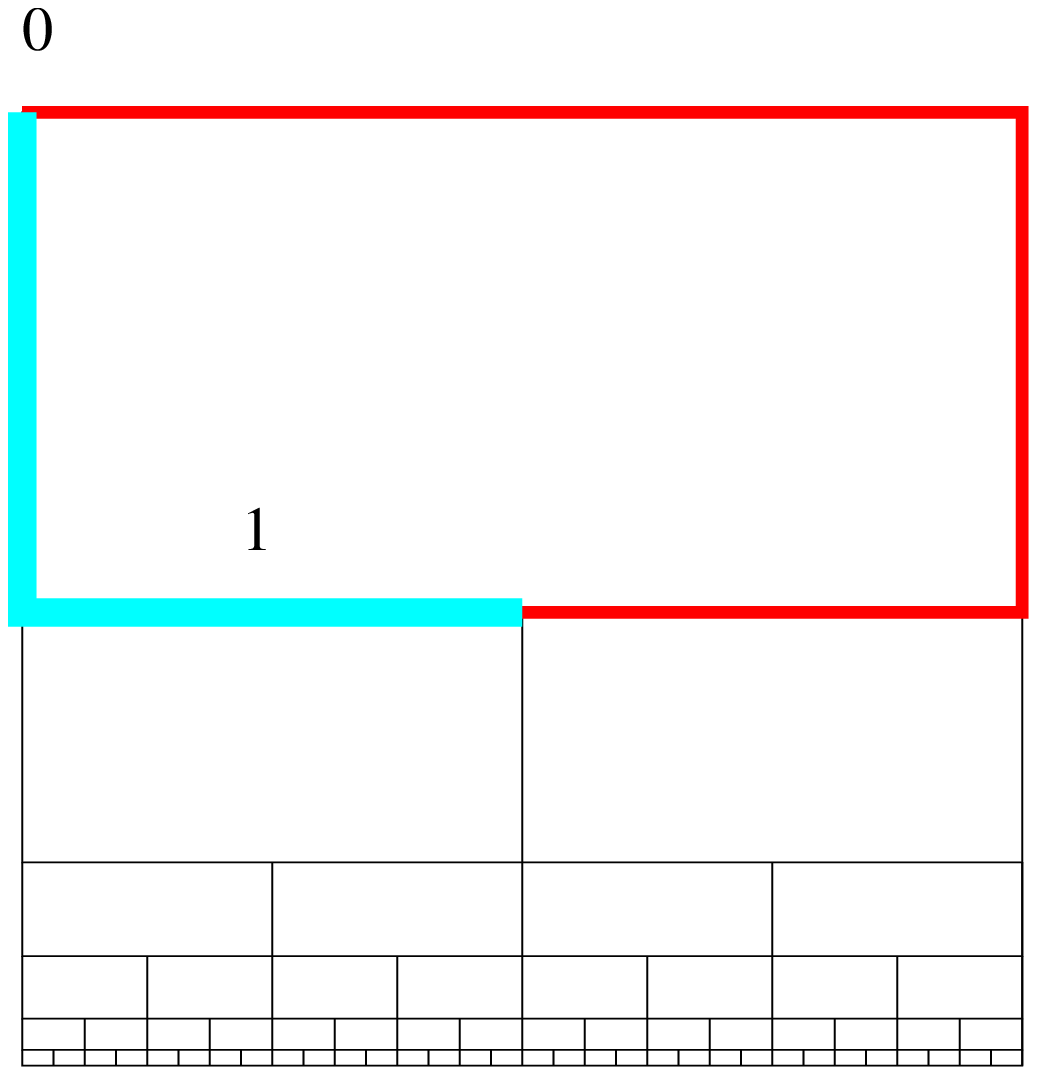}
\et
\et
\ec
  \caption{No $1(-1), (-1)1$ in a run }
  \label{fig:brick}
\end{figure}

\noindent
\textit{Proof.}
For an $N$-run: 
\[
\begin{array}{rclcrlc}
1(-1)& \ra& 01 & \hspace{1cm} & at^{-1}a^{-1} & \ra &t^{-1}a\\ (-1)1
&\ra &0(-1) & \hspace{1cm} & a^{-1}t^{-1}a & \ra& t^{-1}a^{-1}.
\end{array}
\]
For a $P$-run:
\[
\begin{array}{rclcrcl}
1(-1) &\ra &(-1)0 & \hspace{1cm} & ata^{-1}&\ra &a^{-1}t\\ (-1)1 &\ra
&10 & \hspace{1cm} & a^{-1}ta&\ra& at.
\end{array}
\]
See Figure \ref{fig:brick}.
$\Box$

\begin{lem}[No consecutive $11,(-1)(-1)$]\label{lem:No11} 
There exist rewrite rules which do not increase length which can be
applied to a geodesic run to eliminate all occurrences of consecutive
$11$ or $(-1)(-1)$ after the first two entries of an $N$-run and
before the last two entries of a $P$-run.
\end{lem}

\noindent
\textit{Proof.} 
Let $i\in \Z$.\\ For an
$N$-run:
\[
\begin{array}{rclcrlc}
i11& \ra& (i+1)0(-1) & \hspace{1cm} & a^it^{-1}at^{-1}a & \ra &
a^{i+1}t^{-2}a^{-1}\\
 i(-1)(-1) &\ra &(i-1)01 & \hspace{1cm} & a^it^{-1}a^{-1}t^{-1}a^{-1}
& \ra& a^{i-1}t^{-2}a.
\end{array}
\]

\begin{figure}[ht!]
\bc
\bt{ccc}
\bt{c}
\includegraphics[width=6cm]{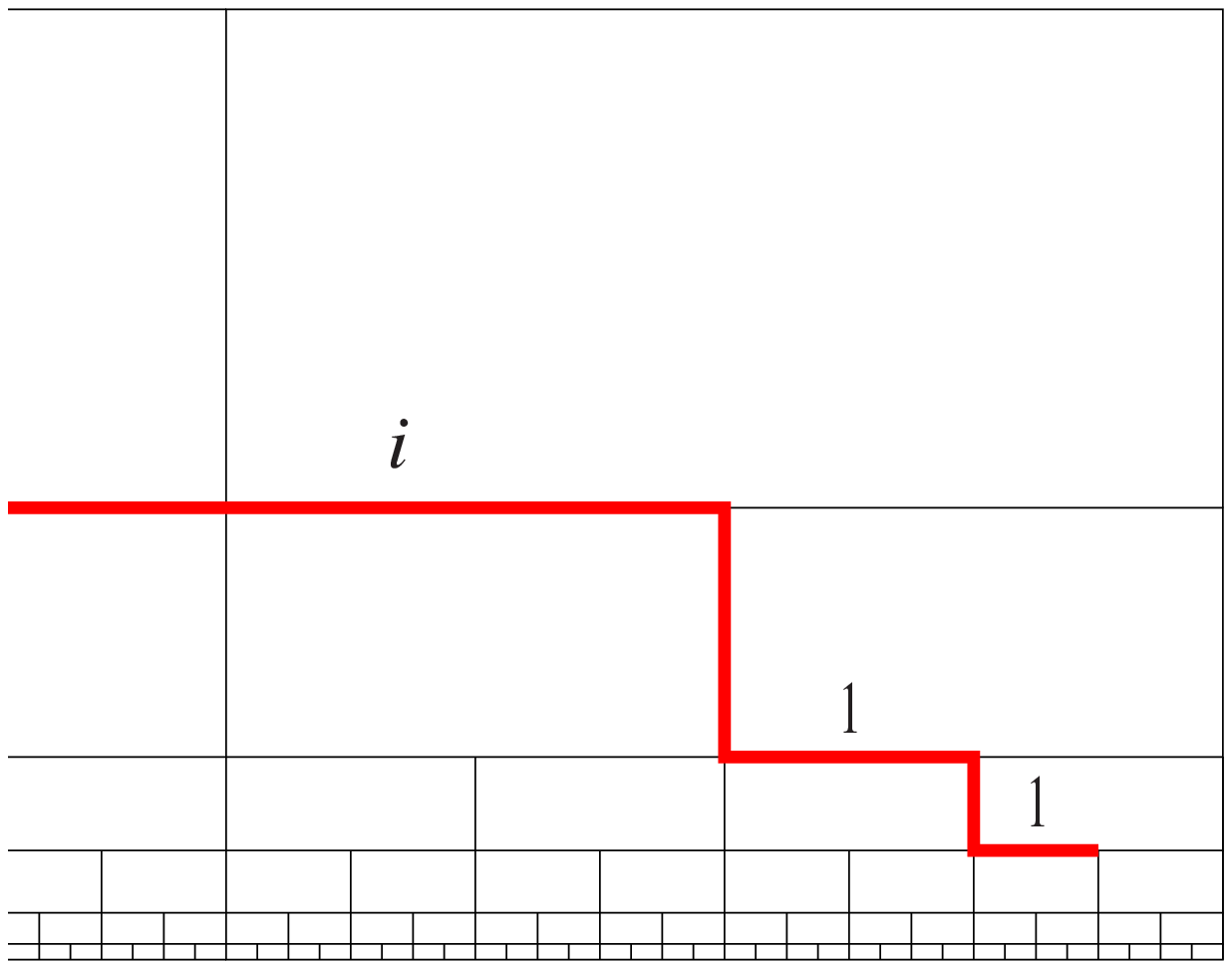}
\et& &
\bt{c}
\includegraphics[width=6cm]{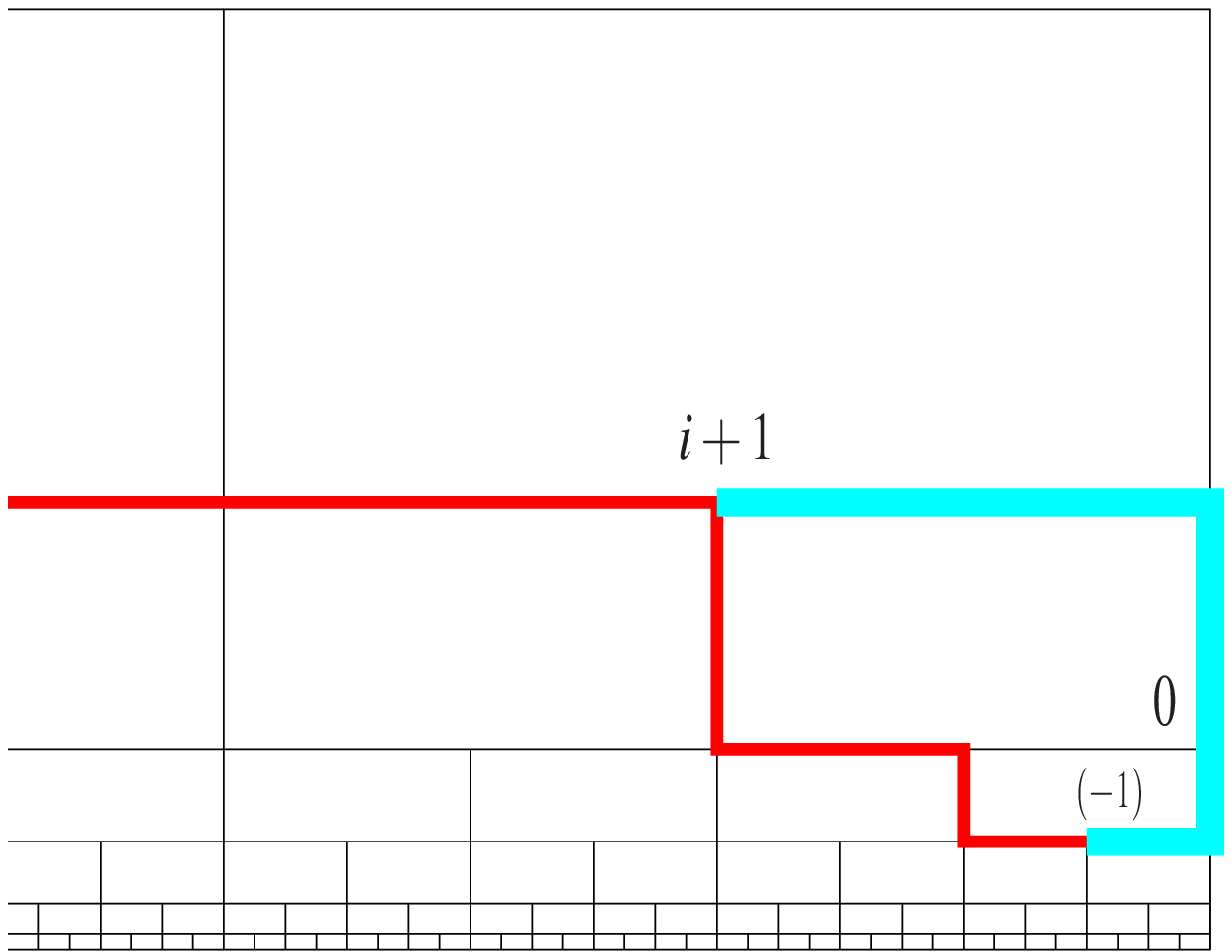}
\et
\et
\ec
  \caption{No $i11, i(-1)(-1)$ in an $N$-run }
  \label{fig:brick1-1}
\end{figure}

These moves are illustrated in Figure \ref{fig:brick1-1}.
 We can always perform these rewrites to get a word of the same length
or shorter. That is, suppose you have an $N$-run, which is geodesic so
we assume has no $1(-1)$ or $(-1)1$. Starting at the right end of the
$N$-run, if there is an $i11$, we know that $i\geq 0$. Replacing this
by $(i+1)0(-1)$ gives a word that is not geodesic if $i>0$, otherwise
gives $10(-1)$. Now if the preceding entry is $(-1)$ the word is not
geodesic, so is $0,1$ or we are at the start of the run. A similar
argument holds when we see $i(-1)(-1)$.

So iterate this procedure until the start of the run is reached. This
eliminates all occurrences of adjacent nonzero entries after the first
two entries.  That is, if the $N$-run starts with $110$ for example,
the rules don't apply.

For a $P$-run:
\[
\begin{array}{rclcrcl}
11i &\ra &(-1)0(i+1) & \hspace{1cm} & atata^i&\ra & a^{-1}t^2a^{i+1}\\
(-1)(-1)i &\ra &10(i-1) & \hspace{1cm} & a^{-1}ta^{-1}ta^i&\ra & at^2a^{i-1}.
\end{array}
\]

Similarly we can always perform these rewrites to get a word of the
same length or shorter, this time starting at the left end of the word
and moving right, so we can eliminate all adjacent nonzero entries
except in the last two positions.  $\Box$

\begin{figure}[ht!]
\bt{ccc}  
\includegraphics[height=3.5cm]{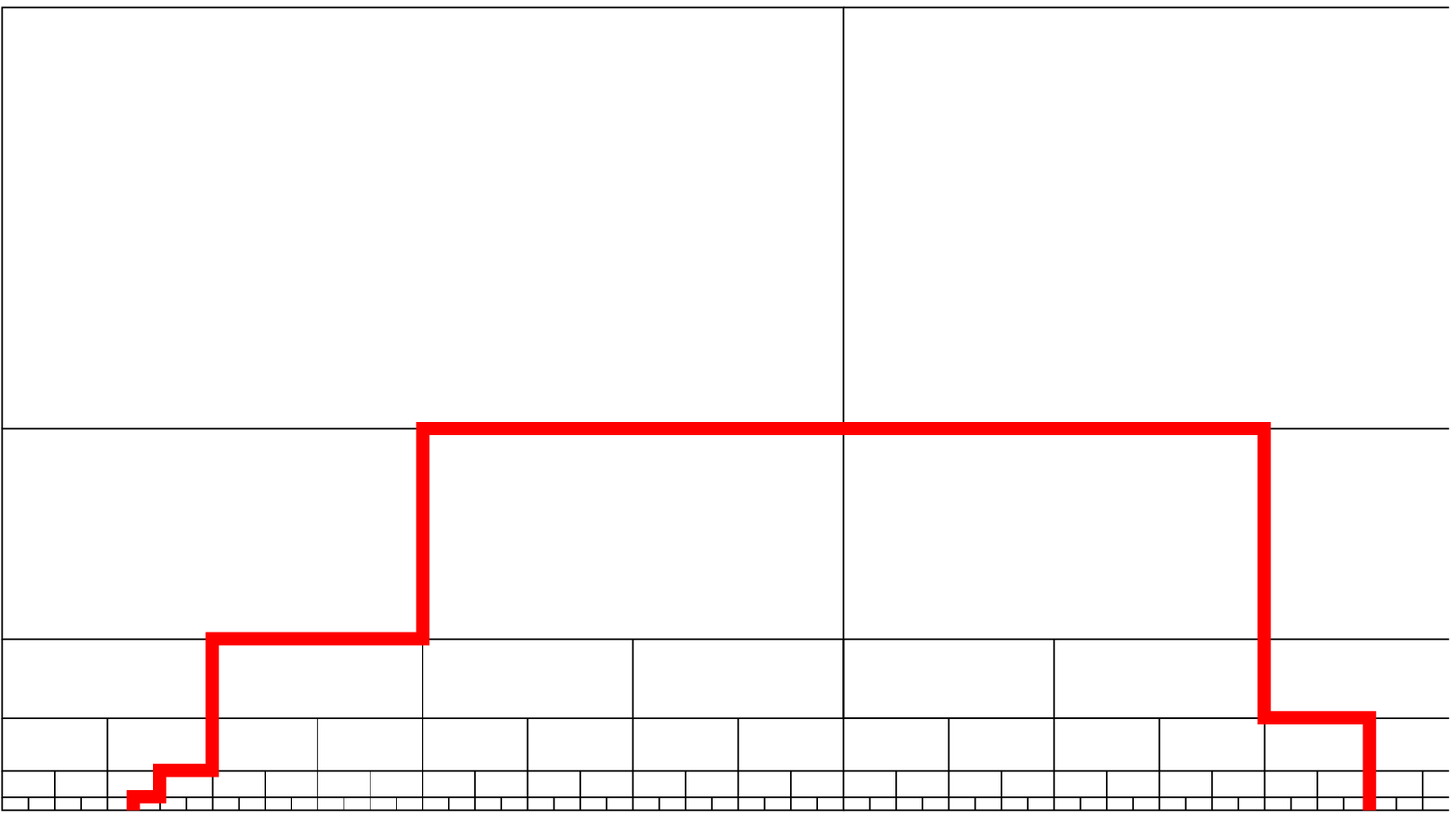}
& &
\includegraphics[height=3.5cm]{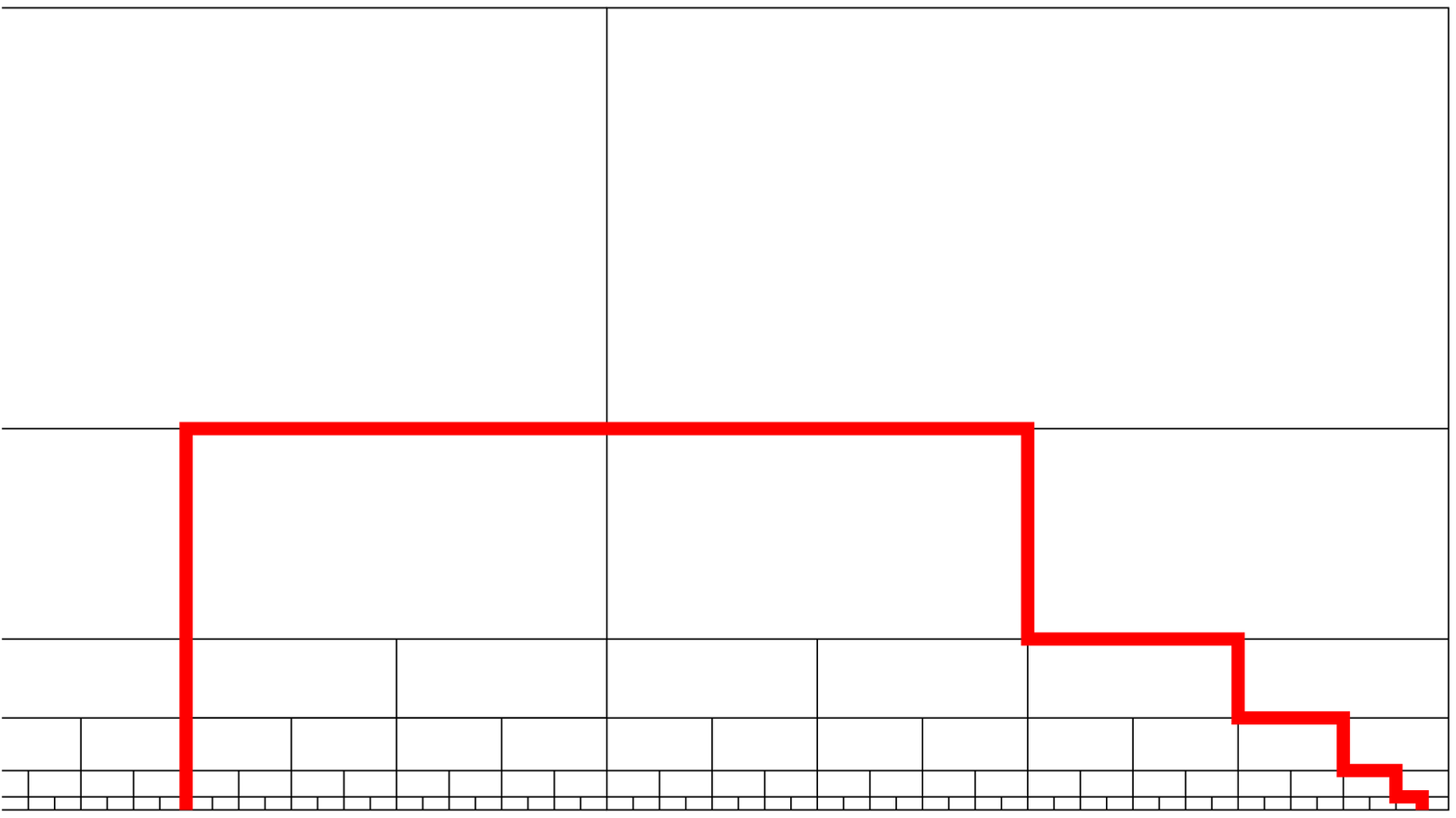}
\et 
  \caption{Pushing an $X$ word to have one non-trivial $N$-run.}
  \label{fig:pushRun}
\end{figure}

Next we will show that every geodesic of one of the ten types can be
``pushed'' into a geodesic word for the same group element that have
at most one non-trivial run. As an example, if
$w=a^{\epsilon_0}ta^{\epsilon_1}t\ldots
a^{\epsilon_k}ta^nt^{-1}a^{\eta_k}t^{-1}\ldots
t^{-1}a^{\eta_1}t^{-1}a^{\eta_0}$ is a geodesic $X$ word, then we can
push the inner subword $a^{\epsilon_k}ta^nt^{-1}a^{\eta_k}$ to
$ta^nt^{-1}a^{\epsilon_k+\eta_k}$, and iteratively push at each level
to get\\ $t^ka^nt^{-1}a^{\epsilon_k+\eta_k}t^{-1}\ldots
t^{-1}a^{\epsilon_1+\eta_1}t^{-1}a^{\epsilon_0+\eta_0}$.
We show this in Figure \ref{fig:pushRun}.

\begin{lem}[At most one run]\label{lem:OneRun}
Every group element is represented by some geodesic of one of the ten
types having at most one non-trivial run.
\end{lem}

\noindent
\textit{Proof.}  By Lemma \ref{lem:TenTypes} each group element is
represented by some geodesic of one of the ten types. If the word is
type $E,N,P$ then there is at most one non-trivial run. If it is $X,
XN$ or $PX$ then by Lemma \ref{lem:commutation} we can push $a$
letters to one side of the $X$ word to get at most one non-trivial
run, as we did in the example above.  For $NP_{\leq}$ words we have
$w=w_Nw_{NP}$ where $w_{NP}$ has zero $t$-exponent, so by Lemma
\ref{lem:commutation} we can push $a$ letters to the left of the $NP$
word to get at most one run non-trivial run. For $XNP$ words we have
$w=w_Xw_Nw_{NP}$ where $w_{NP}$ has zero $t$-exponent, so by Lemma
\ref{lem:commutation} we can push $a$ letters to one side of the $X$
and $NP$ words to get at most one non-trivial run.  For $NP_<$ words
we have $w=w_{NP}w_P$ where $w_{NP}$ has zero $t$-exponent, so by
Lemma \ref{lem:commutation} we can push $a$ letters to the right of
the $NP$ word to get at most one non-trivial run. For $NPX$ words we
have $w=w_Nw_{NP}w_X$ where $w_{NP}$ has zero $t$-exponent, so by
Lemma \ref{lem:commutation} we can push $a$ letters to one side of the
$X$ and $NP$ words to get at most one non-trivial run.  $\Box$

Given that every word can be pushed into a word having at most one
non-trivial run, and we can choose which patterns are not allowed in a
run, we are ready to define the normal form language.

The only issue that remains is the prefix of each run.  For example, a
geodesic of type $X$ can be pushed into a word with exactly one
$N$-run.  The start of this run can be chosen to be either
$a^2t^{-1},a^3t^{-1}, a^{-2}t^{-1}$ or $a^{-3}t^{-1}$, for if the run
starts with $1$ then $tat^{-1}\ra a^2$ so is not geodesic. If it
starts with $4$ or $5$ then by Lemma \ref{lem:TenTypes} $a^4\ra
ta^2t^{-1}$ and $a^5\ra ta^2t^{-1}a$ so we elect to write it starting
with a $2$ instead, and if the run starts with $i\geq 6$ then it is
not geodesic.

  The next few entries could be any one of the following:\\
  $200,201,210,300,301,30(-1),310$ or the negatives of these.

\begin{figure}[ht!]
\bt{ccc}
\includegraphics[height=4cm]{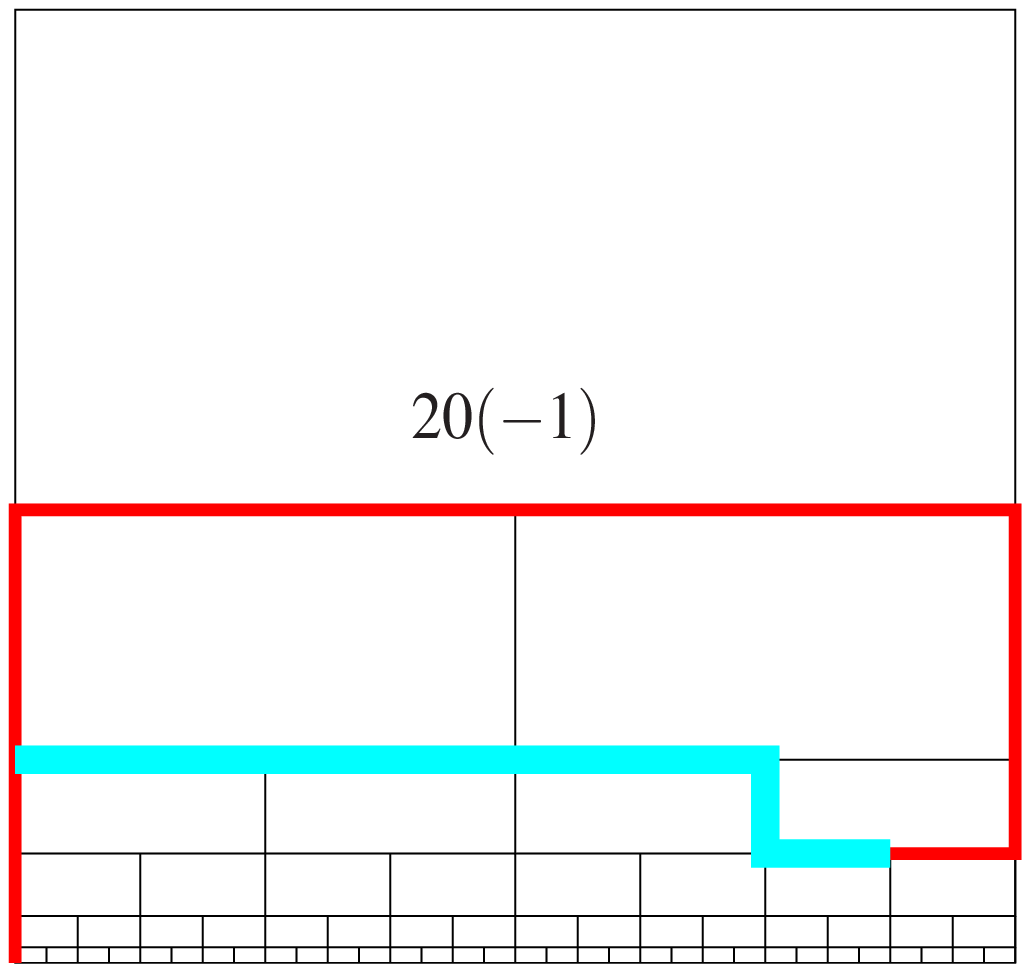}
& &
\includegraphics[height=4cm]{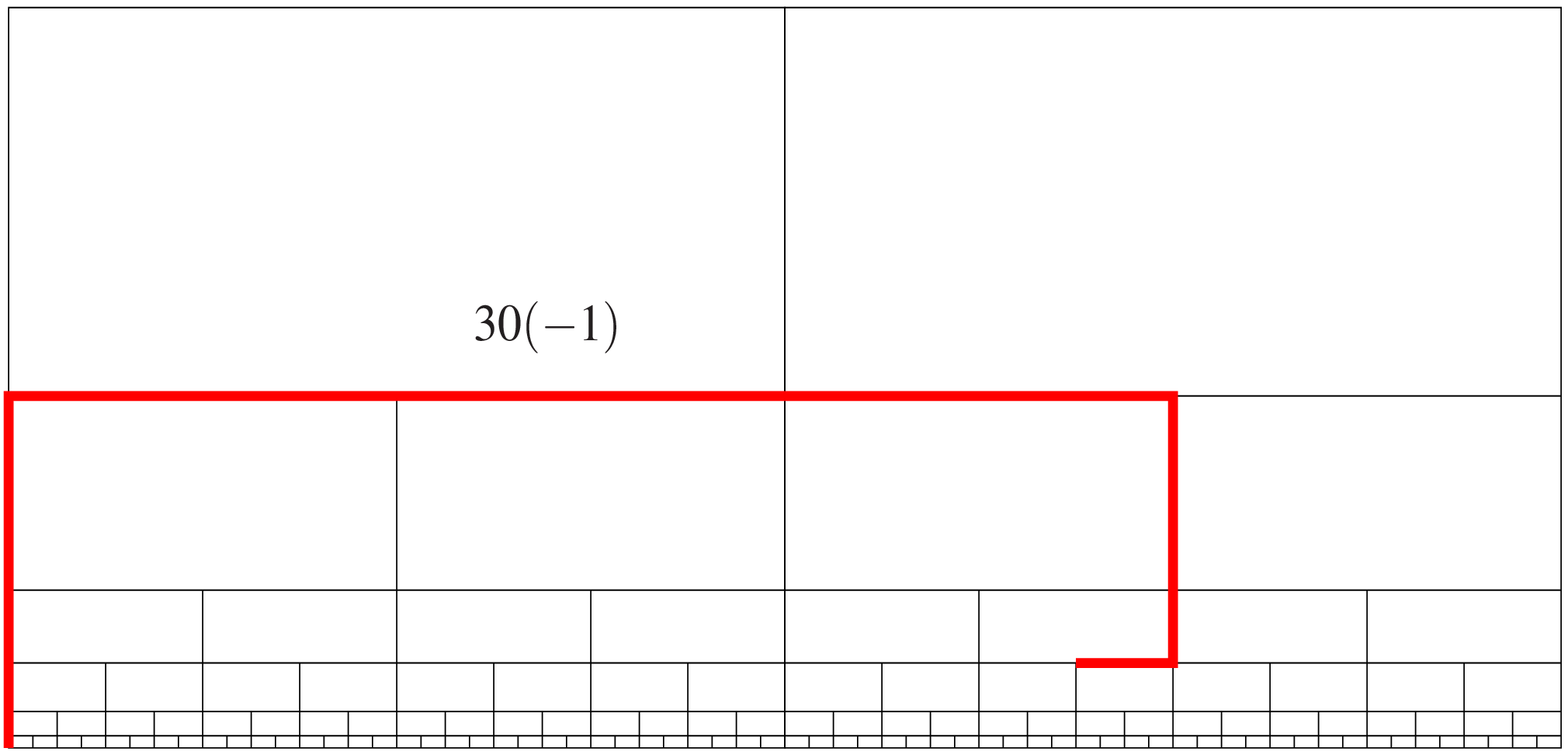}
\et
  \caption{Prefixes for $N$-runs of an $X$ word.}
  \label{fig:prefix}
\end{figure}

Note that the prefix $20(-1)$ is not allowed since $t^2a^2t^{-2}a^{-1}$
is not geodesic, whereas $30(-1)$ is allowed since $t^2a^3t^{-2}a^{-1}$
is geodesic. See Figure \ref{fig:prefix}.

Each case is treated separately in the following lemma. Then after
 these prefixes (suffixes for $P$-runs) the run has only $0, 1,(-1)$
 with no consecutive nonzero entries.

\begin{lem}[Prefixes/suffixes of runs]\label{lem:PrefSuff}
In this lemma we assume that each word has been pushed into a word
with at most one non-trivial run, and that each run has at least three
$t^{\pm 1}$ letters.

\bi
\item The $N$-run in a geodesic word of type $X,XN,XNP$ with
non-positive $t$-exponent sum must start with one of\\
$200,201,210,300,301,30(-1),310$ or the negatives of these.

\item The $N$-run in a geodesic word of type $N,NP_{\leq}$ with
non-positive $t$-exponent sum must start with one of\\
$000,001,010,100,101,10(-1),110,200,201,20(-1),210,300,301,30(-1),310$
or the negatives of these.

\item The $P$-run in a geodesic word of type $P,NP_>$ with positive
$t$-exponent sum must end with one of \\
$000,100,010,001,101,(-1)01,011,002,102,(-1)02,012,003,103,(-1)03,013$
or the negatives of these.

\item The $P$-run in a geodesic word of type $PX,NPX$ with positive
$t$-exponent sum must end with one of\\
$002,102,012,003,103,(-1)03,013$ or the negatives of these.  \ei
\end{lem}

\noindent
\textit{Proof.}  If an $N$-run starts with $i11$ or $i(-1)(-1)$ then
by Lemma \ref{lem:No11} we can replace $i11$ by $(i+1)0(-1)$ and
$i(-1)(-1)$ by $(i-1)01$ without increasing length. Thus the first
three entries of an $N$-run will include a $0$.

If an $N$-run in a word of type $X,XN$ or $XNP$ starts with $i$ with
$|i| \geq 4$ then we can replace $ta^{4+j}t^{-1}$ by
$t^2a^2t^{-1}a^jt^{-1}$ to get a word of the same type and preserving
length.  If an $N$-run in a word of type $X,XN$ or $XNP$ starts with
$i$ with $|i| \leq 1$ then we can replace $ta^it^{-1}$ by $a^{2i}$,
reducing length, contradicting the fact that the word is geodesic.
Thus an $N$-run in a word of type $X,XN$ or $XNP$ starts with
$2,3,(-2)$ or $(-3)$.

This gives the following possibilities for the first three entries:\\
$200,201,20(-1),210,2(-1)0,300,301,30(-1),310,3(-1)0$ or the negatives
of these.  We can eliminate $2(-1)0$ and $3(-1)0$ since they encode
$a^it^{-1}a^{-1}t^{-1}=a^{i-1}t^{-1}at^{-1}$ for $i=2,3$ so are not
geodesic.  We also observe that $20(-1)$ encodes $t^2a^2t^{-2}a^{-1}$
which is not geodesic (as seen in Figure \ref{fig:prefix}).

This leaves $200,201,210,300,301,30(-1),310$ (or their negatives) as
the possible prefixes to the $N$-run in a geodesic of type $X,XN$ or
$XNP$. It is easy to check that each of these prefixes is geodesic.

If the $N$-run in a word of type $N,NP_{\leq}$ starts with $i$ with
$|i| \geq 4$ then we can replace $ta^{4+j}t^{-1}$ by
$t^2a^2t^{-1}a^jt^{-1}$ preserving length. Note that they become words
of type $XN$ or $XNP$.  If the $N$-run in a word of type $N,NP_{\leq}$
starts with $i$ with $|i| \leq 3$ then we can have prefixes of the
form $i0$, $i10$ when $i>0$ and $i(-1)0$ when $i<0$.

Explicitly, this gives\\
$000,001,010,100,101,10(-1),110,200,201,20(-1),210,300,301,30(-1),310$\\
or their negatives. It is easy to check that each of these prefixes is
geodesic. Note that in this case we cannot eliminate $20(-1)$ since
there are no preceding $t$'s.

The proof for $P$-runs follows a similar argument, and is omitted.
$\Box$

\begin{lem}[Short runs]\label{lem:short}
In this lemma we assume that each word has been pushed into a word
with at most one non-trivial run, and that each run has no more than
two $t^{\pm 1}$ letters.
\bi
\item
The geodesics of type $X,XN$ and $XNP$ are the set $L_1$ of words of
the form
\[
\begin{array}{lrl}
ta^it^{-1}a^j, & ij= & (\pm 2)0, (\pm 3)0, 21,31,
(-2)(-1),(-3)(-1); \\ 

ta^it^{-1}a^jt^{-1}a^k, & ijk= & (\pm 2)00, (\pm 3)00, 201,(\pm 3)01,
(-2)0(-1), \\ && (\pm 3)0(-1), 210,310, (-2)(-1)0,(-3)(-1)0; \\

 t^2a^it^{-1}a^jt^{-1}a^k, & ijk= &(\pm 2)00, (\pm 3)00, 201,(\pm
 3)01, (-2)0(-1),  \\ && (\pm 3)0(-1),  210,310, (-2)(-1)0,(-3)(-1)0; \\

ta^it^{-1}a^jt^{-1}a^kt, & ijk= & 201,(\pm 3)01, (-2)0(-1),(\pm
3)0(-1).
\end{array}
\]

\item
The geodesics of type $N$ and $NP_{\leq}$ are the set $L_2$ of words of
the form
\[
\begin{array}{lrl}
a^it^{-1}a^j, & ij= & 00,(\pm 1)0,(\pm 2)0, (\pm 3)0, 0(\pm 1),0(\pm
2), 0(\pm 3),\\ && 11,21,31, (-1)(-1),(-2)(-1),(-3)(-1); \\

a^it^{-1}a^jt, & ij= & 0(\pm 1),0(\pm 2),0(\pm 3),
11,21,31, \\ && (-1)(-1),(-2)(-1),(-3)(-1); \\

a^it^{-1}a^jt^{-1}a^k, & ijk= & 000,(\pm 1)00,(\pm 2)00, (\pm 3)00,
0(\pm 1)0,0(\pm 2)0, 0(\pm 3)0, \\ && 00(\pm 1), 00(\pm 2), 00(\pm 3), (\pm
1)0(\pm 1), (\pm 2)0(\pm 1), (\pm 3)0(\pm 1), \\ &&  110, 210,310,
(-1)(-1)0,(-2)(-1)0,(-3)(-1)0; \\

a^it^{-1}a^jt^{-1}a^kt, & ijk= & 00(\pm 1),00(\pm 2), 00(\pm 3),(\pm
1)0(\pm 1), (\pm 2)0(\pm 1),(\pm 3)0(\pm 1);\\

a^it^{-1}a^jt^{-1}a^kt^2, & ijk= &00(\pm 1),00(\pm 2), 00(\pm 3),(\pm
1)0(\pm 1), (\pm 2)0(\pm 1),(\pm 3)0(\pm 1).
\end{array}
\]

\item
The geodesics of type $P$ and $NP_{>}$ are the set $L_3$ of words of
the form
\[
\begin{array}{lrl}
a^ita^j, & ij= & 00,0(\pm 1),0(\pm 2),0(\pm 3),(\pm 1)0,
11,12,13,\\ && (-1)(-1),(-1)(-2),(-1)(-3);\\

t^{-1}a^ita^j & ij= &(\pm 1)0, 11,12,13,(-1)(-1),(-1)(-2),(-1)(-3);\\

a^ita^jta^k, & ijk= & 000,00(\pm 1),00(\pm 2),00(\pm 3), (\pm 1)0(\pm
1), (\pm 1)0(\pm 2), (\pm 1)0(\pm 3),\\ && 0(\pm
1)0,011,012,013,0(-1)(-1),0(-1)(-2),0(-1)(-3);\\

t^{-1}a^ita^jta^k,  & ijk= &(\pm 1)0(\pm 1),(\pm 1)0(\pm
2),(\pm 1)0(\pm 3);\\

 t^{-2}a^ita^jta^k, & ijk= &(\pm 1)0(\pm 1),(\pm 1)0(\pm 2),(\pm
 1)0(\pm 3).
\end{array}
\]

\item
The geodesics of type $PX$ and $NPX$ (must have positive $t$-exponent)
are the set $L_4$ of words of the form
\[
\begin{array}{lrl}
a^ita^jta^kt^{-1}, & ijk= & 00(\pm 2),00(\pm 3), 012,013, 0(-1)(-2),
0(-1)(-3),\\ &&  102,10(\pm 3), (-1)0(-2),(-1)0(\pm 3).
\end{array}
\]

\ei
\end{lem}

\noindent
\textit{Proof.} The proof is by exhaustive search. For the first two
cases we have either one or two $t^{-1}$ letters, so we consider
$t^pa^it^{-1}a^jt^q$ and $t^pa^it^{-1}a^jt^{-1}a^kt^q$. The
$t$-exponent must be non-positive, so $p+q\leq 1$ in the first case
and $p+q\leq 2$ in the second case. For the $a$-exponents, $|i|\leq 3$
and $|j|,|k|\leq 1$. This gives a finite set of possibilities, so we
run through each and check if it gives a geodesic. Note that the
pattern $20(-1)$ is not a geodesic if it appears in an $N$-run
preceded by a $t$, yet it is geodesic if it is in a $N$ or $NP_{\geq}$
geodesic. 

By Lemma \ref{lem:No11} we choose to reject runs of the form $(i,1,1)$
and $(i,-1,-1)$ in favour of $(i+1,0,-1)$ and $(i-1,0,1)$ respectively,
so that we never see three non-zero entries in a row, even at the
start of a run. The details of the exhaustive check are omitted.

For the third and forth cases we have either one or two $t$ letters,
so we consider $t^{-p}a^ita^jt^{-q}$ and
$t^{-p}a^ita^jta^kt^{-q}$. The $t$-exponent must be positive, so
$p,q=0$ in the third and $p+q\leq 1$ in the forth cases. For the
$a$-exponents, $|k|\leq 3$ and $|j|,|k|\leq 1$. This gives a finite
set of possibilities, so we run through each and check if it gives a
geodesic.

By Lemma \ref{lem:No11} we choose to reject runs of the form $(1,1,i)$
and $(-1,-1,i)$ in favour of $(-1,0,i+1)$ and $(1,0,i-1)$ respectively,
so that we never see three non-zero entries in a row, even at the end
of a run. The details of the exhaustive check are omitted.  $\Box$

\begin{defn}[Normal form]\label{defn:nf} 
There are ten distinct types of \nf\ words.  \bi
\item Type $\NF_E$ words are precisely $\epsilon, a^{\pm 1}, a^{\pm
2}, a^{\pm 3}$.

\item Type $\NF_X, \NF_{XN}$ and $\NF_{XNP}$, all with zero or
negative $t$-exponent, are the words:
$t^ka^{\epsilon_l}t^{-1}a^{\epsilon_{l-1}}t^{-1}\ldots
a^{\epsilon_1}t^{-1}a^{\epsilon_0}t^m$ such that $k>0$ and $l\geq
k+m$, $\epsilon_0 \neq 0$ if $m>0$, the $N$-run starts with one of
$200,201,210,300,301,30(-1),310$  or the negatives of these, and
after this has only $0,1,(-1)$ with no consecutive nonzero entries
(that is, no $1(-1),(-1)1,11$ or $(-1)(-1)$ in the run).

If there are less than three $t^{-1}$ letters in the run, then the
word is in the set $L_1$ of Lemma \ref{lem:short}.

\item Type $\NF_N$ and $ \NF_{NP_{\leq}}$, all with negative
$t$-exponent, are the words:\\
$a^{\epsilon_l}t^{-1}a^{\epsilon_{l-1}}t^{-1}\ldots
a^{\epsilon_1}t^{-1}a^{\epsilon_0}t^k$ such that $0\leq k\leq l$,
$\epsilon_0 \neq 0$ if $k>0$, the $N$-run starts with one of \\
$000,001,010,100,101,10(-1),110,200,201,20(-1),210,300,301,30(-1),310$\\
or the negatives of these, and after this has only $0,1,(-1)$ with no
 consecutive nonzero entries.

If there are less than three $t^{-1}$ letters in the run, then the
word is in the set $L_2$ of Lemma \ref{lem:short}.

\item Type $\NF_P$ and $ \NF_{NP_>}$, all with positive $t$-exponent,
are the words:\\ $t^{-k}a^{\epsilon_0}ta^{\epsilon_{1}}t\ldots
a^{\epsilon_{l-1}}ta^{\epsilon_l}$\\ such that $0\leq k < l$,
$\epsilon_0 \neq 0$ if $k>0$, the $P$-run ends with one of \\
$000,100,010,001,101,(-1)01,011,002,102,(-1)02,012,003,103,(-1)03,013$\\
or the negatives of these, and before this has only $0,1,(-1)$ with no
consecutive nonzero entries.

If there are less than three $t$ letters in the run, then the word is
in the set $L_3$ of Lemma \ref{lem:short}.  

\item Type $\NF_{PX}$ and $ \NF_{NPX}$, all with positive
$t$-exponent, are the words:\\
$t^{-k}a^{\epsilon_0}ta^{\epsilon_{1}}t\ldots
a^{\epsilon_{l-1}}ta^{\epsilon_l}t^{-m}$ such that $k>0, m\geq 0$
and $k+m< l$, $\epsilon_0 \neq 0$ if $k>0$, the $P$-run ends with one
of $002,102,012,003,103,(-1)03,013$ or the negatives of these, and
before this has only $0,1,(-1)$ with no consecutive nonzero entries.

The $P$-run must have at least two $t$ letters since the $t$-exponent
of the word is positive.
If there are less than three $t$ letters in the run, then the word is
in the set $L_4$ of Lemma \ref{lem:short}.  
 \ei
\end{defn}

\begin{lem}[The language of normal forms surjects to the group]
\label{lem:Surject}
Every group element is represented by a normal form word.
\end{lem}

\noindent
\textit{Proof.}  By Lemma \ref{lem:OneRun} every group element is
represented by a geodesic having at most one run. Then by Lemma
\ref{lem:No11} we can remove any occurrences of $11$ and $(-1)(-1)$ in
the run (except possibly at the start of $N$ and $NP_{\leq}$ words and
the end of $P$ and $NP_>$ words) without lengthening the word. Then if
the resulting run does not start (or end) with one of the number
patterns given in Lemma \ref{lem:PrefSuff} relative to its type, it is
not geodesic, and if it does, the word is in normal form.  $\Box$

\begin{defn}[\hnn]
\label{defn:hnn}
If $G$ is a group with presentation $\langle \mc G\;|\; \mc R\rangle$
and $\phi:A\ra B$ is an isomorphism of subgroups $A,B\subseteq G$,
define the {\em \hnn} $G_{\phi}$ of $G$ by $\phi$ to be the group with
presentation $ \langle \mc G, t\;|\; \mc R, \{tat^{-1}=\phi(a) : a \in
A\}\rangle$. The generator $t$ is called the {\em stable letter} and
$A,B$ are called {\em associated subgroups}.
\end{defn}

The group \BS\ is an \hnn\ of $\langle a \rangle$ with the isomorphism
$\phi(a)=a^2$ between associated subgroups $\langle a \rangle$ and
$\langle a^2 \rangle$.  The following fact about \hnn s can be read in
\cite{LS}.  \begin{lem}[Britton's Lemma] \label{lem:Blemma} If $w$ is
a word containing a $t^{\pm 1}$ letter in an \hnn\ of $G_{\phi}$ with
associated subgroups $A,B$ and if $w=_{G_{\phi}}1$ then $w$ must
contain a subword (called a {\em pinch}) of the form $tat^{-1}$ or
$t^{-1}\phi(a)t$ for some element $a\in A$.  \end{lem}

\begin{cor}[$t$-exponent]\label{cor:texp}
For each element $g\in$ \BS\ there is an integer $k$ such that every
word for $g$ has $t$-exponent $k$.
\end{cor}

\noindent \textit{Proof.}  If $w$ represents the identity and has no
$t^{\pm 1}$ letters then its $t$-exponent sum is zero. If $w$
represents the identity and has $t^{\pm 1}$ letters then by
Britton's lemma it contains a pinch. Removing a pinch leaves the
$t$-exponent of $w$ unchanged, so either you can remove all $t^{\pm
1}$ letters, in which case the $t$-exponent sum was zero, or you
cannot remove all $t^{\pm 1}$ letters, in which case the word did not
represent the identity.

If $w$ and $u$ are two words for the same group element with
$t$-exponents $k$ and $l$ respectively, then $wu^{-1}=_{BS} 1$ and has
$t$-exponent $k-l=0$, so $w$ and $u$ have the same $t$-exponent.
$\Box$

\begin{lem}[$a$-exponents]\label{lem:a-exp}
The $X$ word $w=t^ka^jt^{-1}a^{\epsilon_{k-1}}t^{-1}\ldots
t^{-1}a^{\epsilon_0}$ represents the element $a^N$ where
$$N=2^kj+ \sum_{i=0}^{k-1} 2^i\epsilon_i.$$ 

Moreover if each $|\epsilon_i|\leq 1$ for all $i\leq k-1$, $|j|\geq 2$ and
 $\epsilon_{k-1}$ is zero or the same sign as $j$, then $|N|\geq 4$.

Also, the $X$ word $w=a^{\epsilon_0}ta^{\epsilon_1}t\ldots
ta^{\epsilon_{k-1}}ta^jt^{-k}$ represents the element $a^N$ where
$$N=2^kj+ \sum_{i=0}^{k-1} 2^i\epsilon_i,$$ and moreover if each
$|\epsilon_i|\leq 1$ for all $i\leq k-1$, $|j|\geq 2$ and
$\epsilon_{k-1}$ is zero or the same sign as $j$, then $|N|\geq 4$.
\end{lem}

\noindent
\textit{Proof.} 
To prove the first assertion we will use induction on $k$. If $k=1$ we have\\
$w=ta^jt^{-1}a^{\epsilon_0}=a^{2j+\epsilon_0}$.

Assuming the statement holds for $k$, then\\
$w=t^{k+1}a^jt^{-1}a^{\epsilon_k}t^{-1}a^{\epsilon_{k-1}}t^{-1} \ldots
t^{-1}a^{\epsilon_0}$ \\
$=t^ka^{2j+\epsilon_k}t^{-1}a^{\epsilon_{k-1}}t^{-1} \ldots
t^{-1}a^{\epsilon_0}=a^N$ \\ where
$N=2^k(2j+\epsilon_k)+\sum_{i=0}^{k-1} 2^i\epsilon_i$.

The smallest possible value for $|N|$ is when $|j|=2$,
$\epsilon_{k-1}=0$ and each $\epsilon_i$ is $-\frac{j}{|j|}$.  In this
case\\ $|N|\geq 2^k(2)+ 0 + \sum_{i=0}^{k-2} 2^i(-1)$ \\ $= 2^k(2) -
\sum_{i=0}^{k-2} 2^i$\\ $= 2^k(2) - (2^{k-1}-1)$\\ $\geq 2(2)-(1-1)=4$
since $k\geq 1$.

To prove the second assertion we will again use induction on $k$. If
$k=1$ we have\\ $w=a^{\epsilon_0}ta^jt^{-1}=a^{2j+\epsilon_0}$.

Assuming the statement holds for $k$, then\\
$w=a^{\epsilon_0}t\ldots ta^{\epsilon_{k-1}} ta^{\epsilon_k}ta^jt^{-k-1}$\\
$=a^{\epsilon_0}t\ldots ta^{\epsilon_{k-1}} ta^{2^j+\epsilon_k}t^{-k}=a^N$\\
where $N=2^k(2j+\epsilon_k)+\sum_{i=0}^{k-1}2_i\epsilon_i$\\

The smallest possible value for $|N|$ is when $|j|=2$,
$\epsilon_{k-1}=0$ and each $\epsilon_i$ is $-\frac{j}{|j|}$.  In this
case\\ $|N|\geq 2^k(2)+ 0 + \sum_{i=0}^{k-2} 2^i(-1)$ \\ $= 2^k(2) -
\sum_{i=0}^{k-2} 2^i$\\ $= 2^k(2) - (2^{k-1}-1)$\\ $\geq 2(2)-(1-1)=4$
since $k\geq 1$.  $\Box$

%An identical argument gives the same results for
%$w=a^{\epsilon_0}ta^{\epsilon_1}t\ldots
%ta^{\epsilon_{k-1}}ta^jt^{-k}$.  $\Box$

\begin{lem}[Uniqueness for $\NF_E\cup  \NF_X$]
\label{lem:Xunique}
If $w,u\in \NF_E\cup\NF_X$ and $w=_{BS}u$ then $w$ and $u$ are
identical.
\end{lem}

\noindent
\textit{Proof.} If $w,u\in \NF_E$ then $w=a^i$ and $u=a^j$ and
$a^i=_{BS}a^j$ means $a^{i-j}=1$, so $i=j$ and $w$ and $u$ are
identical.

If $w\in \NF_X$ then we can write
$w=t^ka^{\epsilon_k}t^{-1}a^{\epsilon_{k-1}}t^{-1}\ldots
a^{\epsilon_1}t^{-1}a^{\epsilon_0}$ with $k>0$, which evaluates to the
power $N$ with $|N|\geq 4$ by Lemma \ref{lem:a-exp}, so $w$ cannot be
equal to a word in $\NF_E$.

If $u\in \NF_X$ and $w=_{BS}u$ then we can write
$u=t^la^{\eta_l}t^{-1}a^{\eta_{l-1}}t^{-1}\ldots
t^{-1}a^{\eta_{k}}t^{-1}\ldots a^{\eta_1}t^{-1}a^{\eta_0}$, where
without loss of generality we are assuming that $k\leq l$.  Since both
words evaluate to the same power of $a$ we have
\begin{eqnarray*}
\epsilon_k2^k+\epsilon_{k-1}2^{k-1}+\ldots +\epsilon_12+\epsilon_0 & =
&\eta_l2^l+\eta_{l-1}2^{l-1}+\ldots+\eta_{k}2^{k}+\ldots
+\eta_12+\eta_0.\end{eqnarray*}

Let $i\in \mathbb N$ such that $\epsilon_j=\eta_j$ for all $j<i$ and
$\epsilon_i\neq \eta_i$. Then cancelling and dividing through by $2^i$
we have
\begin{equation}\label{eqn1}
\epsilon_k2^{k-i}+\epsilon_{k-1}2^{k-1-i}+\ldots +\epsilon_i =
\eta_l2^{l-i}+\eta_{l-1}2^{l-1-i}+\ldots +\eta_i.
\end{equation}

%If $i=k$ then $\epsilon_i=2$ or $3$ and since there is only one \nf\
%word for $a^2$ and $a^3$ then $w$ and $u$ are identical.

If $i=k$ then $|\epsilon_k|=2$ or $3$ and we have
$\epsilon_k=\eta_l2^{l-k}+\eta_{l-1}2^{l-1-k}+\ldots +\eta_k$.  If
$l=k$ then $\epsilon_k=\eta_k$ so $w$ and $u$ are identical.  If
$l\geq k+1$ then
$|\epsilon_k|=|\eta_l2^{l-k}+\eta_{l-1}2^{l-k-1}+\ldots \eta_k|\geq
4$ since $|\eta_l|\geq 2$ and $\eta_{l-1}$ is either $0$ or the same
sign as $\eta_l$, but $|\epsilon_k|\leq 3$ so this is a contradiction.

If $i<k$ then $\epsilon_i, \eta_i$ are either $0,\pm 1$ since they
occur in the middle of a run. By Equation \ref{eqn1} they must be of
the same parity, and they cannot both be zero so one is $1$ and one is
$(-1)$. If $i+1<k$ then $\epsilon_{i+1}=\eta_{i+1}=0$ and we
contradict the equation since one side is equal to $1 \mod 4$ and the
other is $(-1) \mod 4$.

So $i+1=k$, so the run in $w$ starts with $210$ or $310$ (or their
negatives).  Then $w=t^ka^st^{-1}aw''$ and $u=t^ku't^{-1}a^{-1}w''$
with $s=2,3$ so $u'=_{BS}a^{s+1}$ so is $a^3$ or $a^4$, which by
Lemma \ref{lem:TenTypes} is written as $ta^2t^{-1}$ if it occurs
in a normal form word. Then the run in $u$ must start with either
$3(-1)$ or $20(-1)$, neither of which is allowed in a \nf\ word, so
$w$ and $u$ are identical.  $\Box$

\begin{lem}[Uniqueness for $\NF_N\cup \NF_{XN}$]
\label{lem:Nunique}
If $w,u\in \NF_N\cup\NF_{XN}$ and $w=_{BS}u$ then $w$ and $u$ are
identical.
\end{lem}

\noindent
\textit{Proof.}  If $w$ and $u$ are two \nf\ words representing the
same group element, then they have the same $t$-exponent by Lemma
\ref{cor:texp}.  If $w,u\in \NF_{XN}$ with $t$-exponent $(-k)$ then
$t^kw,t^ku$ are in $\NF_{X}$ so by Lemma \ref{lem:Xunique} they are
identical. Note that $\NF_{XN}$ and $\NF_X$ words have the same
$N$-run structure, the only difference is the length of the $t^l$
prefix.

If $w\in \NF_N$ then let $w=a^{\epsilon_k}t^{-1}\ldots
t^{-1}a^{\epsilon_0}$ and let \\ $u=u't^{-1}a^{\eta_{k-1}}t^{-1}\ldots
t^{-1}a^{\eta_0}$ where $u'$ evaluates to $a^n$ and is type $X$ or
$E$.  The words $t^kw$ and $t^ku$ evaluate to the same power of $a$,
which is $\epsilon_k2^k+\ldots
+\epsilon_0=n2^k+\eta_{k-1}2^{k-1}+\ldots +\eta_0$. Let $i\in \mathbb
N$ such that $\epsilon_j=\eta_j$ for all $j<i$ and $\epsilon_i\neq
\eta_i$.  Then cancelling and dividing through by $2^i$ we get
\begin{equation}\label{eqn2}
\epsilon_k2^{k-i}+\ldots
+\epsilon_i=n2^{k-i}+\eta_{k-1}2^{k-i-1}+\ldots +\eta_i.
\end{equation}

If $i=k$ then $\epsilon_k=n$.  Now $|\epsilon_k|\leq 3$ and $u'$ is an
$E$ or $X$ word with the same $a$-exponent.  By Lemma \ref{lem:a-exp}
if $u'$ is type $X$ then it evaluates to $a^N$ with $|N|\geq 4$, so
$u'$ is type $E$, indeed it is exactly $a^{\epsilon_k}$, so $w$ and
$u$ are identical.

If $i<k$ then $\epsilon_i,\eta_i=\pm 1$ since they are in the middle
of a run, and have the same parity by Equation \ref{eqn2}. If $i<k+1$
then we have a contradiction since $\epsilon_{i+1}=\eta_{i+1}=0$ and
the equation has $4x+1$ on one side and $4y-1$ on the other for
integers $x,y$.  So $i=k+1$ and
$\epsilon_k2+\epsilon_{k-1}=n2+\eta_{k-1}$ so $n=\epsilon_k\pm 1$
since $\epsilon_{k-1} -\eta_{k-1}=\pm 2$, and $\epsilon_{k-1}$ has the
same sign as $\epsilon_k$.

If $u'$ is type $X$ then $|n|\geq 4$ by Lemma \ref{lem:a-exp} but
$|\epsilon_k|\leq 3$, so the only chance for equality is when the run
in $w$ starts with $31$ and $\eta_{k-1}=-1$. Then $u'=_{BS}a^4$ which
is written as $ta^2t^{-1}$ in a normal form word, but then the run in
$u$ starts with $20(-1)$ which is not allowed.  Thus $u$ is also in
$\NF_N$. Without loss of generality assume $\epsilon_k>0$ so
$\epsilon_{k-1}=1$ and $\eta_{k-1}=-1$. Then $n$ must be negative
since the run in $u$ starts with $n(-1)$, and we have a contradiction.
$\Box$

\begin{lem}[Uniqueness for $\NF_P\cup \NF_{PX}$]
\label{lem:Punique}
If $w,u\in \NF_P\cup\NF_{PX}$ and $w=_{BS}u$ then $w$ and $u$ are
identical.
\end{lem}

\noindent
\textit{Proof.}  If $w,u\in \NF_P\cup\NF_{PX}$ then $w^{-1}$ and
$u^{-1}$ are in $\NF_N\cup\NF_{XN}$, so by Lemma \ref{lem:Nunique}
since $w^{-1}=_{BS}u^{-1}$ then $w^{-1}$ and $u^{-1}$ are identical,
and so $w$ and $u$ are identical.  $\Box$

\begin{lem}[Uniqueness]\label{lem:Unique}
Every group element is represented by a unique normal form word.
\end{lem}

\noindent
\textit{Proof.}  If $w$ and $u$ are two \nf\ words representing the
same group element, then they have the same $t$-exponent by Lemma
\ref{cor:texp}.

 If $w$ and $u$ have zero $t$-exponent then they are of the form
$E,X,NP$ or $XNP$.  If neither is $NP$ or $XNP$ then they are
identical by Lemma \ref{lem:Xunique}.  If one is $NP$ or $XNP$ then
let $w=w't^{-1}a^{\epsilon_{k-1}}t^{-1}\ldots t^{-1}a^{\epsilon_0}t^k$
and $u=u't^{-1}a^{\eta_{l-1}}t^{-1}\ldots t^{-1}a^{\eta_0}t^l$ where
$w',u'$ evaluate to powers of $a$ and assume without loss of
generality that $k>0$ and $k\geq l$.  Then
$wu^{-1}=w't^{-1}a^{\epsilon_{k-1}}t^{-1}\ldots
t^{-1}a^{\epsilon_0}t^{k-l}a^{-\eta_0}t\ldots
ta^{-\eta_{l-1}}(u')^{-1}$\\ $=_{BS}1$. Since $k>0$ then
$\epsilon_0=\pm 1$ so if we replace $w'$ and $u'$ by the corresponding
powers of $a$ (by pinching $ta^st^{-1}$ subwords) we have a word that
does not admit any pinches, contradicting Britton's Lemma. Thus
$k=l$. Then the words $wt^{-k}$ and $ut^{-k}$ are equal and in
$\NF_N\cup\NF_{XN}$ so by Lemma \ref{lem:Nunique} must be identical,
so $w$ and $u$ are identical.

 If $w$ and $u$ have negative $t$-exponent then they are of the form
$N,XN,NP$ or $XNP$.  If neither is $NP$ or $XNP$ then they are
identical by Lemma \ref{lem:Nunique}.  If one is $NP$ or $XNP$ then
let $w=w't^{-1}a^{\epsilon_{k-1}}t^{-1}\ldots t^{-1}a^{\epsilon_0}t^l$
and let \\ $u=u't^{-1}a^{\eta_{p-1}}t^{-1}\ldots t^{-1}a^{\eta_0}t^q$
where $k>l,p>q$, and $w',u'$ evaluate to powers of $a$. Assume without
loss of generality that $l>0$ and $l\geq q$.  Then\\
$wu^{-1}=w't^{-1}a^{\epsilon_{k-1}}t^{-1}\ldots
t^{-1}a^{\epsilon_0}t^{l-q}a^{-\eta_0}t\ldots
ta^{-\eta_{p-1}}(u')^{-1}=_{BS}1$. Since $l>0$ then $\epsilon_0=\pm 1$
so after replacing $w'$ and $u'$ by the corresponding powers of $a$,
we have a word that does not admit any more pinches, contradicting
Britton's Lemma. Thus $l=q$. Then the words $wt^{-l}$ and $ut^{-l}$
are equal and in $\NF_N\cup\NF_{XN}$ so by Lemma \ref{lem:Nunique}
must be identical, so $w$ and $u$ are identical.

 If $w$ and $u$ have positive $t$-exponent then they are of the form
$P,PX,NP$ or $NPX$.  If neither is $NP$ or $NPX$ then they are
identical by Lemma \ref{lem:Punique}.  If one is $NP$ or $NPX$ then
let $w=t^{-k}a^{\epsilon_0}t\ldots t a^{\epsilon_{k-1}}w'$ and let
$u=t^{-p}a^{\eta_0}t \ldots ta^{\eta_{k-1}tu'}$ where $k<l,p<q$, and
$w',u'$ evaluate to powers of $a$. Assume without loss of generality
that $k>0$ and $k\geq p$.  Then \\$u^{-1}w=
(u')^{-1}t^{-1}a^{-\eta_{k-1}}t^{-1}\ldots t^{-1}a^{-\eta_0}
t^{p-k}a^{\epsilon_0}t\ldots t a^{\epsilon_{k-1}}w'=_{BS}1$.  Since
$k>0$ then $\epsilon_0=\pm 1$ so after replacing $w'$ and $u'$ by
their corresponding powers of $a$ we have a word that cannot be
pinched, contradicting Britton's Lemma. Thus $k=p$. Then the words
$t^kw$ and $t^ku$ are equal and in $\NF_P\cup\NF_{PX}$ so by Lemma
\ref{lem:Punique} must be identical, so $w$ and $u$ are identical.
$\Box$

\begin{lem}[Normal forms are geodesic]\label{lem:geodNF} 
Each \nf\ word is a geodesic.
\end{lem}

\noindent
\textit{Proof.}  Suppose that a word $w\in \NF$ is not
geodesic. Choose a geodesic word $u=_{BS}w$ that is one of the ten
types in Lemma \ref{lem:TenTypes}. By Lemma \ref{lem:OneRun} we can
move $u$ into a word $u'$ of the same length having one run.

If $u'$ is in \nf\ then since $w$ and $u'$ are both \nf\ words that
equate to the same group element then $w,u'$ must be identical by Lemma
\ref{lem:Unique}.

If $u'$ is not in \nf, it either violates the prefix rules (as in
Lemma \ref{lem:PrefSuff}) or has an adjacent pair of nonzero digits in
its run.  

If the run in $u'$ has an occurrence of $1(-1)$ or $(-1)1$ then $u'$ is
not geodesic.  If the run in $u'$ has an occurrence of $11$ or
$(-1)(-1)$ that is not at the start of an $N$-run or the end of a
$P$-run, then by Lemma \ref{lem:No11} we can perform a length
preserving rewrite to eliminate it. If this causes $u'$ to have a
$1(-1)$ then $u$ was not geodesic, and it it causes $u'$ to have a
$11$ or $(-1)(-1)$ then repeatedly applying Lemma \ref{lem:No11} from
right to left in an $N$-run, or left to right in a $P$-run, we can
eliminate all occurrences of pairs of nonzero digits.

Finally if the start or end is not one of the prefixes in Lemma
\ref{lem:PrefSuff} then either $u'$ is not geodesic (if the prefix is
$20(-1)$ for example), or is equal to a normal form word of the same
length, which means that the original word $w$ is geodesic.  $\Box$

\section{The main theorem}\label{sect:mainthm}

\begin{thm}
The language $\NF$ is a 1-counter language.
\end{thm}

\noindent
\textit{Proof.} The ten types of normal-form geodesics listed in
Definition \ref{defn:nf} break up into five cases.  The set $\NF_E$ is
a 1-counter \lan\ since it is finite. We can describe a $\Z$-automaton
for each of the remaining four cases to accept the remaining nine
types.

Consider the set of normal forms words of type $X,XN$ and $XNP$. The
language $L_1$ of Lemma \ref{lem:short} describes the set of normal
form words of these types with at most two $t^{-1}$ letters in the
$N$-run, and since $L_1$ is finite, it is a regular language.

Let $L_1'$ be the set of words of the form
$\{t^ka^it^{-2}at^{-1},t^ka^jt^{-2}a^{-1}t^{-1}\: | \; k=1,2,3, i=2, \pm
3, j=-2,\pm 3\}$. This is a finite set so is regular, and is the set
of $X$ (and $XN$) normal form words with three $t^{-1}$'s in the
$N$-run, that corresponds to the prefix $201,301,30(-1) $ and their
negatives.

\begin{figure}[ht!]
\bt{ccc}
\includegraphics[height=6.5cm]{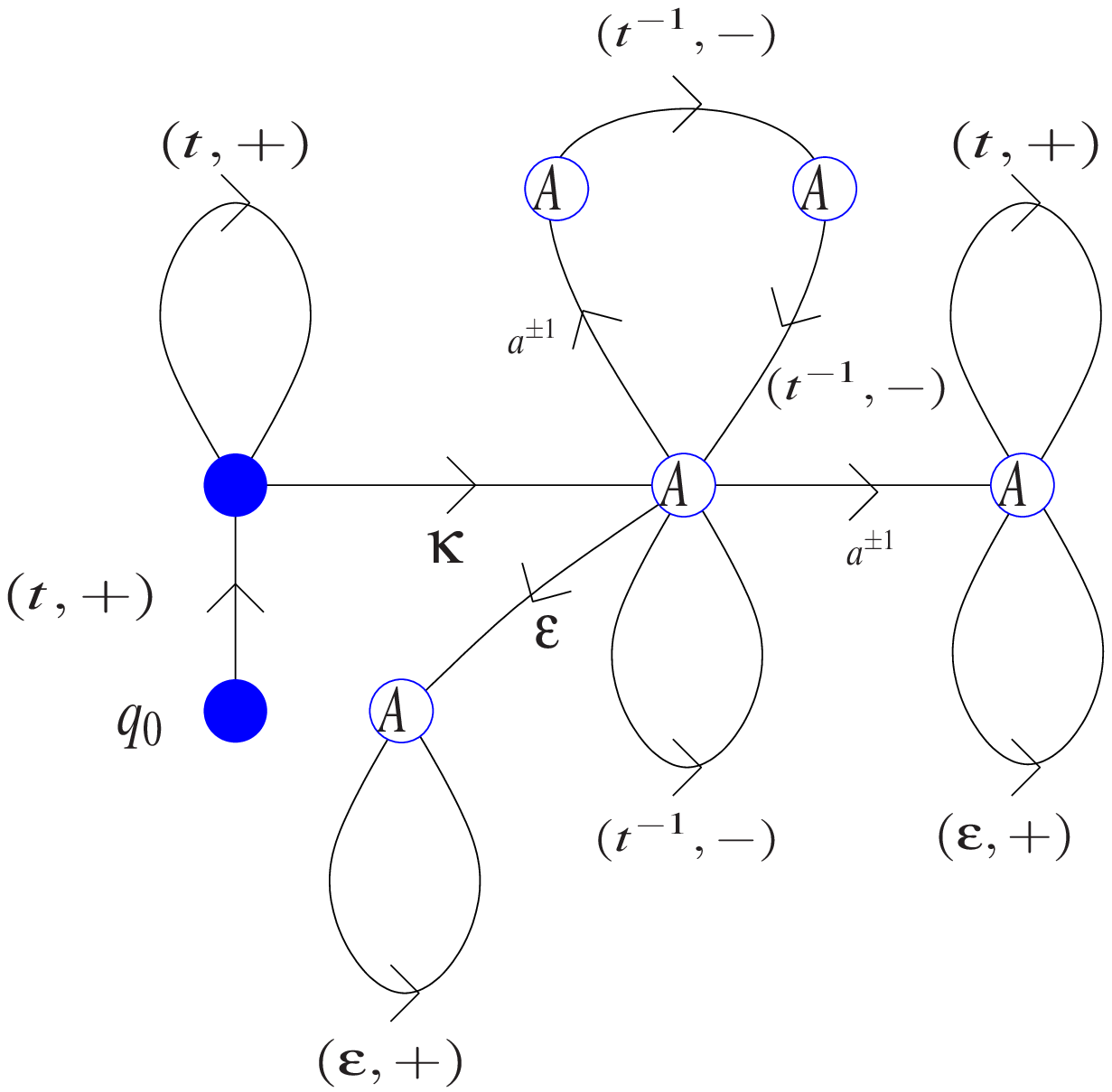}& & 
\includegraphics[height=6.5cm]{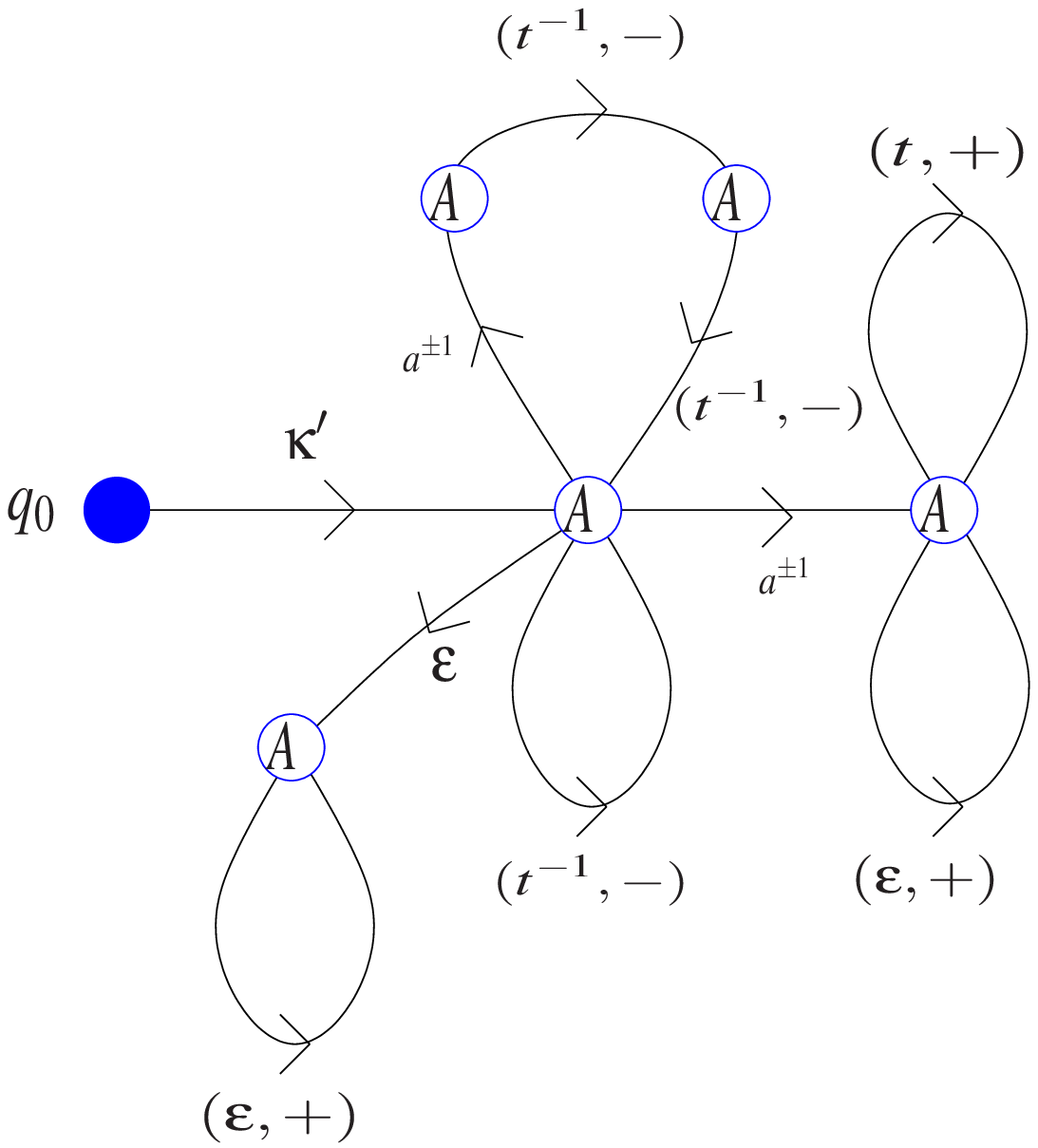}
\et
  \caption{Counter automata for normal form $X,XN,XNP$ words and
  $N,NP_{\leq}$ words with $N$-run length at least $3$.}
  \label{fig:aut1}
\end{figure}

The remaining $X,XN$ and $XNP$ normal form words (with an $N$-run of
$3$ or more $t^{-1}$ letters) are accepted by the automaton on the
left of Figure \ref{fig:aut1}.  The edge labeled $\kappa$ stands for a
collection of paths labeled by
\[
\begin{array}{ll}
 a^i(t^{-1},-)(t^{-1},-)(t^{-1},-), & i=\pm 2,\pm 3; \\
 a^i(t^{-1},-)(t^{-1},-)a(t^{-1},-)(t^{-1},-),   & i= 2, \pm 3;\\
 a^i(t^{-1},-)(t^{-1},-)a^{-1}(t^{-1},-)(t^{-1},-),   & i= -2, \pm 3;\\
 a^i (t^{-1},-)a(t^{-1},-)(t^{-1},-), & i= 2, 3; \\
 a^i (t^{-1},-)a^{-1}(t^{-1},-)(t^{-1},-), & i= -2, -3. 
\end{array}
\]
The union of these three (regular and 1-counter) languages is 1-counter.

Next, consider the set of normal forms words of type $N$ and
$NP_{\leq}$.  The language $L_2$ of Lemma \ref{lem:short} describes
the set of normal form words of these types with at most two $t^{-1}$
letters in the $N$-run, and since $L_2$ is finite, it is a regular
language.

Let $L_2'$ be the set of words of the form $\{a^it^{-2}a^{\pm
1}t^{-1}\: | \;  i=0, \pm 1,\pm 2, \pm 3\}$. This is a
finite set so is regular, and is the set of $N$ (and $NP_{\leq}$) normal
form words with three $t^{-1}$'s in the $N$-run, that corresponds to
the prefix $00(\pm 1),10(\pm 1),20(\pm 1),30(\pm 1) $ and their
negatives.

The remaining $N$ and $NP_{\leq}$ normal form words (with an $N$-run
of $3$ or more $t^{-1}$ letters) are accepted by the automaton on the
right of Figure \ref{fig:aut1}.  The edge labeled $\kappa'$ stands for
a collection of paths labeled by
\[
\begin{array}{ll}
 a^i(t^{-1},-)(t^{-1},-)(t^{-1},-), & i=0, \pm 1, \pm 2,\pm 3; \\
 a^i(t^{-1},-)(t^{-1},-)at^{-1}(t^{-1},-),   & i= 0,\pm 1,\pm 2, \pm 3;\\
 a^i(t^{-1},-)(t^{-1},-)a^{-1}(t^{-1},-)(t^{-1},-),& i= 0,\pm 1,\pm 2, \pm 3;\\
 a^i (t^{-1},-)a(t^{-1},-)(t^{-1},-), & i= 0,1,2, 3; \\
 a^i (t^{-1},-)a^{-1}(t^{-1},-)(t^{-1},-), & i= 0,-1,-2, -3. 
\end{array}
\]

Next, consider the set of normal forms words of type $P$ and $NP_>$. The
language $L_3$ of Lemma \ref{lem:short} describes the set of normal
form words of these types with at most two $t$ letters in the
$P$-run, and since $L_3$ is finite, it is a regular language.

Let $L_3'$ be the set of words of the form $\{ta^{\pm 1}t^2a^i \: | \;
i=0, \pm 1,\pm 2, \pm 3\}$.  This is a finite set so is regular, and
is the set of $P$ (and $NP_>$) normal form words with three $t$'s in
the $P$-run, that corresponds to the suffix $(\pm 1)00,(\pm 1)01,(\pm
1)02,(\pm 1)03 $ and their negatives.

\begin{figure}[ht!]
\bt{ccc}
\includegraphics[height=5cm]{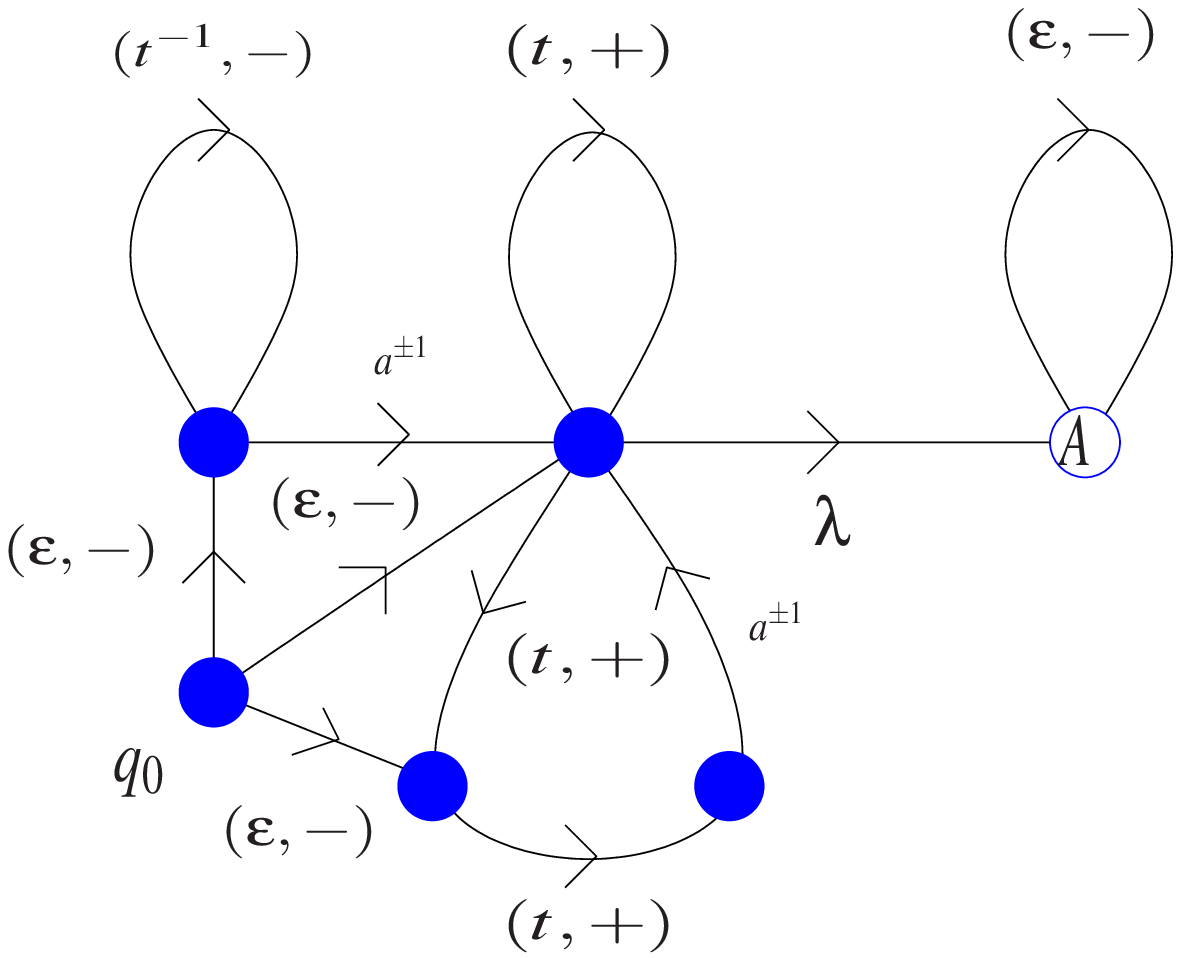}& & 
\includegraphics[height=5cm]{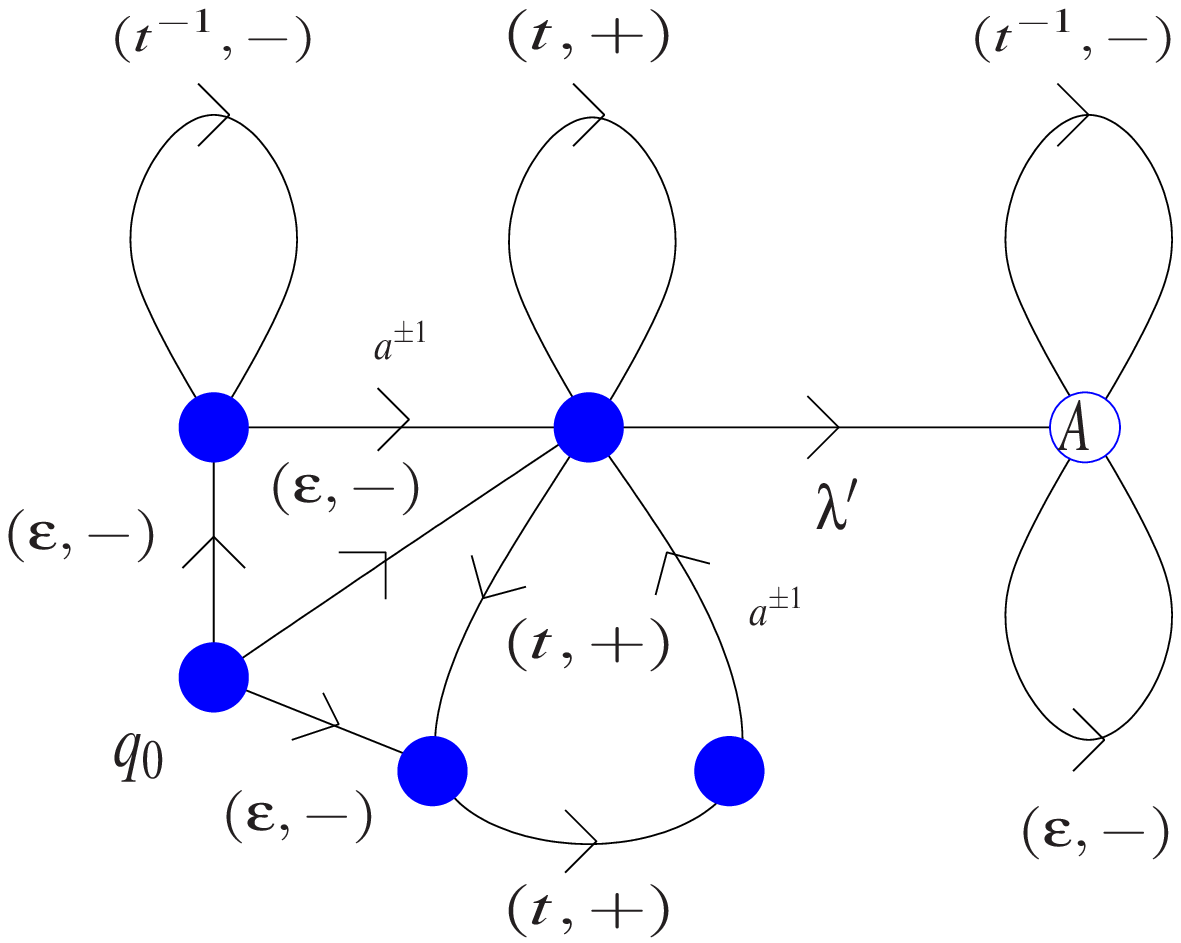}
\et
  \caption{Counter automata for normal form $P,NP_>$ words and
  $PX,NPX$ words with $P$-run length at least $3$.}
  \label{fig:POS}
\end{figure}

The remaining $P$ and $NP_>$ normal form words (with a $P$-run of $3$
or more $t$ letters) are accepted by the automaton on the left of
Figure \ref{fig:POS}.  The edge labeled $\lambda$ stands for a
collection of paths labeled by
\[
\begin{array}{ll}
(t,+)(t,+)(t,+)a^i & i=0, \pm 1, \pm 2,\pm 3; \\ 
(t,+)(t,+)a(t,+)(t,+)a^i & i= 0,\pm 1,\pm 2, \pm 3;\\ 
(t,+)(t,+)a^{-1}(t,+)(t,+)a^i & i= 0,\pm 1,\pm 2, \pm 3;\\ 
(t,+)(t,+)a(t,+)a^i & i= 0,1,2, 3;\\ 
(t,+)(t,+)a^{-1}(t,+)a^i & i= 0,-1,-2, -3.
\end{array}
\]

Lastly, consider the set of normal forms words of type $PX$ and
$NPX$.  The language $L_4$ of Lemma \ref{lem:short} describes
the set of normal form words of these types with (at most) two $t$
letters in the $P$-run, and since $L_4$ is finite, it is a regular
language.

Let $L_4'$ be the set of words of the form
$\{tat^2a^it^{-k},ta^{-1}t^2a^jt^{-k},
\: | \; k=1,2,3, i=2, \pm
3, j=-2,\pm 3\}$. This is a finite set so is regular, and is the set
of $PX$ (and $NPX$) normal form words with three $t$'s in the
$P$-run, that corresponds to the suffix $102,103,(-1)03 $ and their
negatives.

The remaining $PX$ and $NPX$ normal form words (with a $P$-run of $3$ or more
$t$ letters) are accepted by the automaton on the right of Figure
\ref{fig:POS}.  The edge labeled $\lambda'$ stands for a collection of
paths labeled by
\[
\begin{array}{ll}
(t,+)(t,+)(t,+) a^i & i=\pm 2,\pm 3; \\
 (t,+)(t,+)a(t,+)(t,+)a^i   & i= 2, \pm 3;\\
 (t,+)(t,+)a^{-1}(t,+)(t,+)a^i   & i= -2, \pm 3;\\
 (t,+)(t,+)a(t,+)a^i  & i= 2, 3; \\
 (t,+)(t,+)a^{-1}(t,+)a^i & i= -2, -3. 
\end{array}
\]

By Lemma \ref{lem:closure1counter} the union of a 1-counter and
a regular language is 1-counter so each of the ten types is 1-counter,
and by Lemma \ref{cor:union} the union of 1-counter languages is
1-counter.  $\Box$

\begin{cor}
The language of normal forms for BS$(1,2)$ with the standard
generating set is context-free.
\end{cor}

%\section{Remarks}
%With a more careful normal form description one can prove the same
%results for BS$(1,p)$ when $p\geq 3$. The ideas in Lemmas
%\ref{lem:commutation}, \ref{lem:No-11} and \ref{lem:No11}
%can be applied to understanding geodesics for words that evaluate to
%powers of $a$ in the non-solvable case as well.

%Same results for BS$(1,p)$ with $p\geq 3$. Normal form description is 
%more involved. The idea of pushing $a$ letters left and right through
% $X$ words applies in the non-solvable case as well.
% so some kind of rules on what patterns are allowed in $X$ words might 
%prove useful in finding the growth functions for non-solvable BS groups.

\section{Full language of geodesics}

In this section we prove that the language of all geodesic words in
the standard generating set is not counter.  To prove this we will
mimic the proof of Theorem \ref{thm:cfnotcounter}. Recall that in
that proof we constructed a word $ww^R$ on three symbols whose prefix
is square-free and suffix is its reverse, and applied the Swapping
Lemma (Lemma \ref{lem:swap}) to obtain a contradiction.

Let $w$ be a word in \BS\ with no $a^{-1}$ letters.  Define the {\em
$t$-encoding} of $w$ to be a string of integers $n_1n_2\ldots n_k$
such that $w=t^{n_1}at^{n_2}\ldots at^{n_k}$. If $w$ starts (or
respectively ends) with an $a$ then $n_1=0$ (or respectively $n_k=0$).

As an example, the word 
\begin{eqnarray*}
 at^2a^2ta^3t^4at^{-9}at^2at^{-1} & = & t^0 a t^2 a t^0 a t a t^0
a t^0 a t^4 a t^{-9} a t^2 a t^{-1}
\end{eqnarray*}
is encoded as $0201004(-9)2(-1).$ Note that previously our encodings
have been of $a$-exponents, but this new encoding will be useful for
the argument to follow.

\begin{figure}[ht!]
\bc
\bt{c}
\includegraphics[height=4cm]{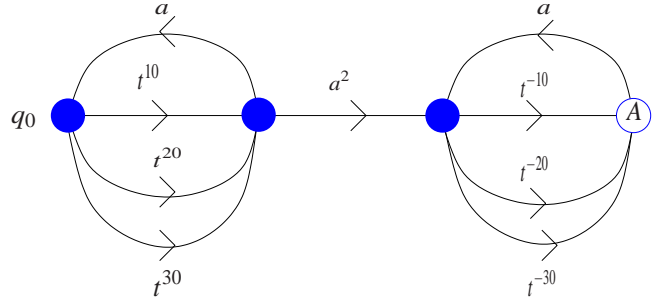}
\et
\ec
  \caption{A finite state automaton accepting the language $L$ in the
  proof of Theorem \ref{thm:BSFullnotcounter}}
  \label{fig:FSA_tencodings}
\end{figure}

\begin{thm}\label{thm:BSFullnotcounter}
The language of all geodesic words in \BS\ with respect to the
generating set $\{a^{\pm 1}, t^{\pm 1}\}$ is not counter.
\end{thm}

\noindent
\textit{Proof.}  Suppose that the full language is counter, and call
it $C$.  Define $L$ to be the set of words in $\{a,t^{\pm 1}\}$
accepted by the \fsa\ in Figure \ref{fig:FSA_tencodings}.  That is,
$L$ is the set of $PN$ words whose $t$-encodings are words of the form
$\{10,20,30\}\{10,20,30\}^*0\{-10,-20,-30\}\{-10,-20,-30\}^*.$

Since $L$ is regular, the intersection of $C$ and $L$ is counter.  Let
$M$ be a counter automaton accepting $C\cap L$, with alphabet $a^{\pm
1},t^{\pm 1}$.  We can construct a new counter automaton $M'$ which
accepts the set of $t$-encoded words of $C\cap L$ as follows.

The states, start state, accept states and counters are the same as
for $M$. The new alphabet is $\{0,\pm 10,\pm 20,\pm 30\}$.
The transitions are defined as follows.

If there is a path labelled by $t^ia$ in $M$ from $p$ to $q$, then add
an edge in $M'$ from $p$ to $q$ labeled by $i$, and the counters are
changed by the same amount as they were following the path $t^ia$ in
$M$. Thus a word is accepted by $M$ if and only if its encoding is
accepted by $M'$. Since $M$ accepts $C\cap L$, the only subwords of
the form $t^ia$ that appear in accepted words are for $i=0,\pm 10,\pm
20$ or $\pm 30$.  Let $p$ be the swapping length for $M'$.

 Next, take a Thue-Morse word in three symbols, which we choose to be
$10,20,30$, of length greater than $2p$. This word encodes a $P$ word
$u$ of some $t$-exponent $10c$.  We wish to find some kind of
``reverse'' of $u$, as we did in the proof of Theorem
\ref{thm:cfnotcounter}.  We find a word $v$ to act as the ``reverse''
by the following procedure.

\be
\item Write $u$ as $t^{10}a^{\epsilon_1}t^{10}a^{\epsilon_2}\ldots
t^{10}a^{\epsilon_k}t^{10}$ where $\epsilon_i=0,1$.
\item Reverse this word.
\item Replace $a^0$ with $a^1$ and $a^1$ with $a^0$ in this word.
\item Replace $t^{10}$ with $t^{-10}$ in this word to get $v$.  \ee

For example, the Thue-Morse word $10,20,30,10,30,20,10,20,30,20,10,30$
encodes the word
\begin{eqnarray*}
u & = & t^{10}at^{20}at^{30}at^{10}at^{30}at^{20}at^{10}at^{20}a
t^{30}at^{20}at^{10}at^{30}
\end{eqnarray*}

\noindent
Step 1: Write $u$ as
\begin{eqnarray*}
u & = &
%t^{10}a^{1}t^{10}a^0t^{10}a^{1}t^{10}a^0t^{10}a^0t^{10}a^{1}t^{10}a^{1}t^{10}
%a^0t^{10}a^0t^{10}a^{1}t^{10}a^0t^{10}
%a^{1}t^{10}a^{1}t^{10}a^0t^{10}a^{1}t^{10}a^0t^{10}a^0t^{10}a^{1}t^{10}a^0
%t^{10}a^{1}t^{10}a^{1}t^{10}a^0t^{10}a^0t^{10}
%\end{eqnarray*}
|a^{1}|a^0|a^{1}|a^0|a^0|a^{1}|a^{1}|a^0|a^0|a^{1}|a^0|
a^{1}|a^{1}|a^0|a^{1}|a^0|a^0|a^{1}|a^0|a^{1}|a^{1}|a^0|a^0|
\end{eqnarray*}
where the $t^{10}$ terms are replaced by bars $|$, to make it easier
to read.

\noindent
Step 2: Reversing this word gives 
\begin{eqnarray*}
u^R & = &
|a^0|a^0|a^1|a^1|a^0|a^1|a^0|a^0|a^1|a^0|a^1|a^1|a^0|a^1|a^0|a^0|
a^1|a^1|a^0|a^0|a^1|a^0|a^1|.
\end{eqnarray*}

\noindent Step 3: Replacing $a^0$ by $a^1$ and vice versa gives
\begin{eqnarray*}
 & & |a^1|a^1|a^0|a^0|a^1|a^0|a^1|a^1|a^0|a^1|a^0|a^0|a^1|a^0|a^1|a^1|
a^0|a^0|a^1|a^1|a^0|a^1|a^0|.
\end{eqnarray*}

\noindent
Step 4: Replacing $t^{10}$ by $t^{-10}$ gives 
\begin{eqnarray*}
v & = & \dag a^1 \dag a^1 \dag a^0 \dag a^0 \dag a^1 \dag a^0 \dag a^1
\dag a^1 \dag a^0 \dag a^1 \dag a^0 \dag a^0 \dag a^1 \dag a^0 \dag
a^1 \dag a^1 \dag a^0 \dag a^0 \dag a^1 \dag a^1 \dag a^0 \dag a^1
\dag a^0 {\dag}\\ & = & t^{-10}at^{-10}at^{-30} at^{-20}at^{-10}a
t^{-20}a t^{-30}at^{-20}a t^{-10}at^{-30}a t^{-10}at^{-20}at^{-20}
\end{eqnarray*}

\noindent
where $\dag$ represents $t^{-10}$.

The $t$-encoding for $v$ is then 
\begin{equation*}
(-10)(-10)(-30)(-20)(-10)(-20)(-30) (-20)(-10)(-30)(-10)(-20)(-20).
\end{equation*}
 Note that $v$ does not have to be square-free. Note also that the
 $t$-exponent of $v$ is $-10c$, where $10c$ is the $t$-exponent of
 $u$.

\begin{figure}[ht!]
\bt{c}
\includegraphics[width=14.5cm]{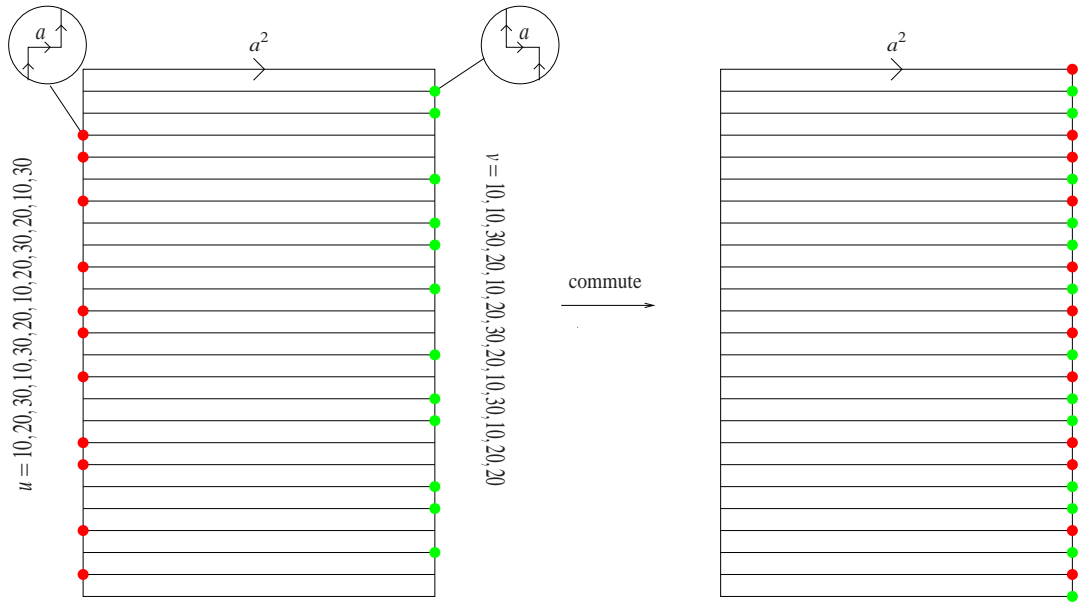}
\et
  \caption{The word $w=ua^2v$ drawn in a sheet of the \cg.}
  \label{fig:sheet-tencodings}
\end{figure}

Now to understand what motivated us to produce this $v$ from $u$,
consider the word $w=ua^2v=uat^0av$. This word is type $X$.  Drawing
$w$ in a sheet of the \cg\ we see that at every tenth level there is
an $a$ letter, either on the part going up the sheet (the $u$ part) or
the part going down (the $v$ part). See the left side of Figure
\ref{fig:sheet-tencodings}.

We will now show that $w$ is a geodesic. Consider the word $w'$
obtained from $w$ by commuting all $a$ letters to the right. Since
there is exactly one $a$ at every tenth level of $w$, we have
$w'=t^{10c}a^2t^{-10}(at^{-10})^{c-1}.$ Then $w'$ is a normal form $X$
word, since its $N$-run is of the form $200\ldots$ with no consecutive
non-zero entries.  Thus by Lemma \ref{lem:geodNF} is geodesic, and
since $w'$ has the same length as $w$ then $w$ is geodesic. So $w$ is
in $C\cap L$, it is accepted by the counter automaton $M$, and its
$t$-encoding is accepted by $M'$.

Applying the Swapping Lemma (Lemma \ref{lem:swap}) to the encoding of
$w$, we switch two adjacent subwords in the first half of $w$, that
is, in the $t$-encoding of $u$, which is square-free.

\begin{figure}[ht!]
\bt{c}
\includegraphics[width=14.5cm]{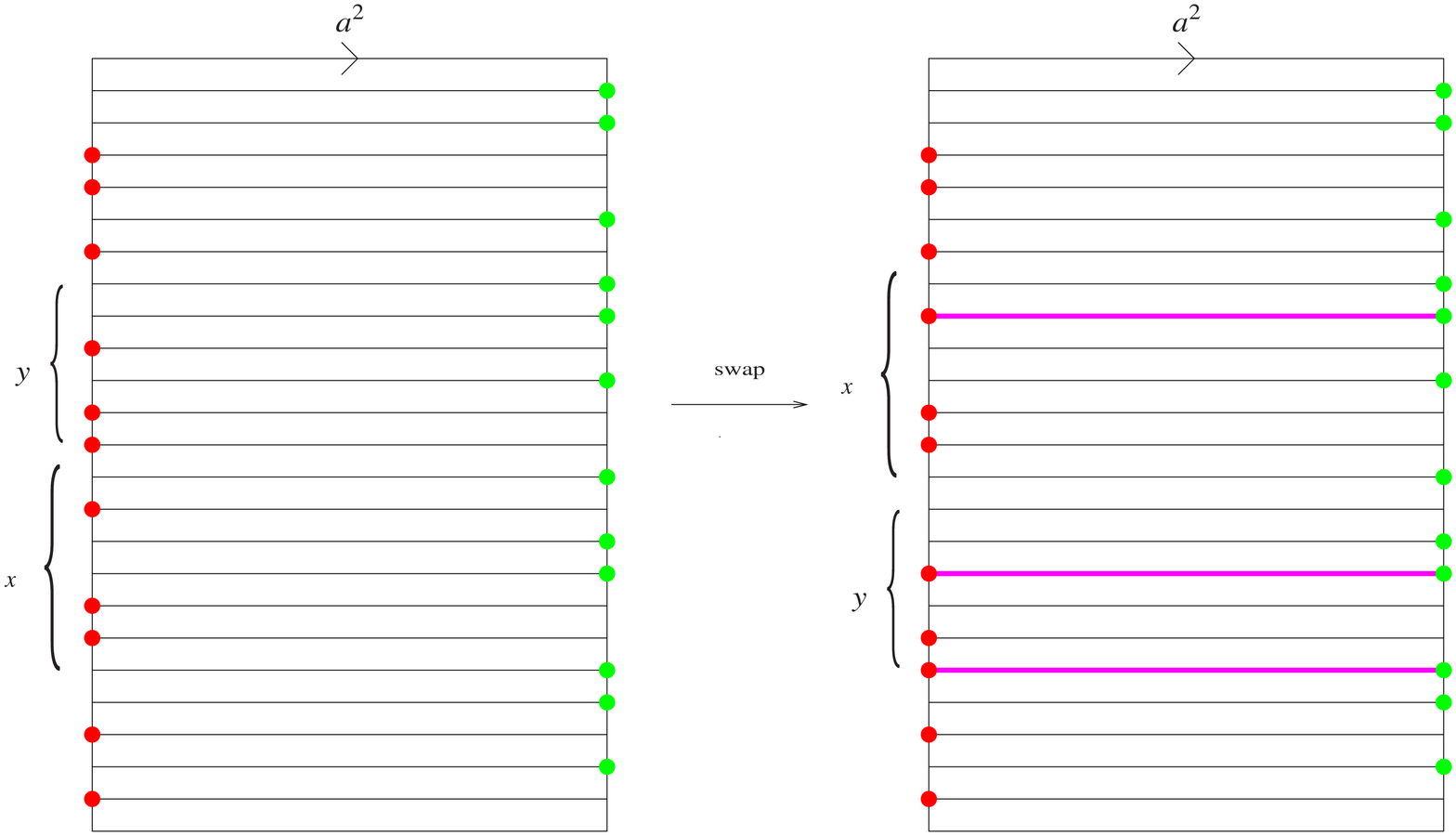}
\et
  \caption{Swapping two subwords in the $P$ part of $w$ leads to a
  word with $t^{-1}a^2t^{-1}$.}
  \label{fig:sheet-swap}
\end{figure}

This new string is a $t$-encoding of some other word in the group,
which is an $X$ word, essentially the same as $w$ except that at some
level(s) we see a shift one step to the right in both sides of the
word (viewed in the sheet of the \cg). See Figure
\ref{fig:sheet-swap}.

When we commute $a$-letters to the right in this word, we will see
$t^{-1}a^2t^{-1}$ at some point(s) in the $N$-run, and thus the swapped
word is not a geodesic, so not in $C\cap L$, and this is a
contradiction.  $\Box$

%Question - is the full language of geodesics context-free?

\section{Acknowledgements}

My sincere thanks to Bob Gilman, Ray Cho, Walter Neumann, Jon
McCammond, Susan Hermiller, Sarah Rees, Rick Thomas, Nik Ruskuc, Kim
Ruane, Mauricio Gutierrez, Sean Cleary, Jennifer Taback and Gretchen
Ostheimer for their help and suggestions that have all contributed to
this work. I wish to thank the reviewer of this paper for pointing out
that the normal form language described here is a 1-counter language,
as well as many other very useful suggestions and corrections.  The
labels for the figures were produced using Andrew Rechnitzer's
\texttt{equation$\_$edit} program.

\end{document}